\documentclass[conference]{IEEEtran}

\newcommand{\bol}{\boldsymbol}

\newcommand{\ney}{\boldsymbol{y}}                          
\newcommand{\nex}{\boldsymbol{x}}

\newcommand{\nez}{\boldsymbol{z}}

\newcommand{\ner}{\mathbf{r}}                      
\newcommand{\TE}{\mathrm{TE}}                           
\newcommand{\TM}{\mathrm{TM}}
\newcommand{\de}{\,\mathrm{d}}                               
\newcommand{\e}{\operatorname{e}}                               
                     
\newcommand{\inc}{\mathrm{inc}}

\newcommand{\andtext}{\quad\mbox{and}\quad}

\newcommand{\p}{\partial}

\newcommand{\real}{\mathrm{Re}\,}    
                                   
\newcommand{\imag}{\mathrm{Im}\,}

\newcommand{\lf}{\left}
\newcommand{\rg}{\right}

\newcommand{\R}{\mathbb{R}}

\newcommand{\mgf}{\bold H}                                        
\newcommand{\elf}{\bold E}

\newcommand{\nor}{\hat{\bold n}}

\newcommand{\bracenom}{\genfrac{\lbrace}{\rbrace}{0pt}{}}

\usepackage{multirow}
\usepackage{amsmath}
\usepackage{amsfonts}
\usepackage{amsmath}
\usepackage{amssymb}  
\usepackage{graphicx}
\usepackage{caption}
\usepackage{mathrsfs}
\usepackage{upgreek}
\usepackage{amsthm}
\usepackage{subfig}
\usepackage{booktabs}
\usepackage{authblk}

\usepackage[noadjust]{cite}
\usepackage[T1]{fontenc}

\usepackage[dvipsnames]{xcolor}
\usepackage[normalem]{ulem}

\usepackage{cite}
\usepackage{url}

\ifCLASSINFOpdf
\else
\fi
\usepackage{authblk}

\title{Windowed Green Function MoM for Second-Kind Surface Integral Equation Formulations of Layered Media Electromagnetic Scattering Problems}
\author[1]{Rodrigo Arrieta}
\author[2]{Carlos~P\'erez-Arancibia}
\affil[1]{\small Department of Electrical Engineering, PUC Chile ({\tt riarrieta@uc.cl})}
\affil[2]{\small Department of Applied Mathematics, University of Twente ({\tt c.a.perezarancibia@utwente.nl})}

\begin{document}

\maketitle
\thispagestyle{plain}
\pagestyle{plain}

\begin{abstract}
 This paper presents a second-kind surface integral equation method for the numerical solution of frequency-domain electromagnetic scattering problems by locally perturbed layered media in three spatial dimensions. Unlike standard approaches, the proposed methodology does not involve the use of layer Green functions. It instead leverages an indirect M\"uller formulation in terms of free-space Green functions that entails integration over the entire unbounded penetrable boundary. The integral equation domain is effectively reduced to a small-area surface by means of the windowed Green function method, which exhibits high-order convergence as the size of the truncated surface increases. The resulting (second-kind) windowed integral equation is then numerically solved by means of the standard Galerkin method of moments (MoM) using RWG basis functions. The methodology is validated by comparison with Mie-series and Sommerfeld-integral exact solutions as well as against a layer Green function-based MoM. Challenging examples including realistic structures relevant to the design of plasmonic solar cells and all-dielectric metasurfaces, demonstrate the applicability, efficiency, and accuracy of the proposed methodology. 

\end{abstract}


\IEEEpeerreviewmaketitle

\section{Introduction}\label{sec:intro}
\IEEEPARstart{P}{roblems} of electromagnetic (EM) scattering and radiation in the presence of planar layered media have played an important role in the development of electromagnetics since the beginning of the 20th century, when the seminal works of Zenneck and Sommerfeld on the propagation of radio waves over the surface of the earth appeared~\cite{Michalski:2016kk}. Their relevance lies in that in many application areas it is crucial to determine the scattering from localized perturbations (e.g., surface roughness, small-size inclusions, meta-atoms) and/or the field produced by localized sources (e.g., antenna feeds) embedded within physically large structures that, away from a certain region of interest, can be effectively assumed as planar and infinite (e.g., the surface of the earth, silicon substrates). This is often the case in numerous problems in radio communications~\cite{tamir1967radio}, remote subsurface sensing~\cite{song2005reconstruction}, microwave circuits~\cite{ling1999efficient,mosig1985general}, nano-optical metamaterials~\cite{yu2014flat},  photonics~\cite{saleh2019fundamentals}, and plasmonics~\cite{maier2005plasmonics}.

Popular numerical approaches to layered media scattering include differential equation-based methods, such as the finite difference~\cite{taflove1995computational} and finite element methods~\cite{jin2015finite}, and surface integral equation (SIE) methods, such as the method of moments (MoM)~\cite{Harrington1993} (also known as the boundary element method) and Nystr\"om methods~\cite{tong2020nystrom}.  Unlike differential equation-based methods, SIE methods do not suffer from dispersion errors, and they can easily handle unbounded domains and radiation conditions at infinity without recourse to perfectly matched layers  or approximate absorbing/transparent boundary conditions for truncation of the computational domain. Additionally, SIE methods rely on discretization of the relevant physical boundaries, and they therefore give rise to linear systems of reduced dimensionality which, although dense, can be efficiently solved by means of iterative solvers in conjunction with fast algorithms~\cite{chew2001fast}. 

\begin{figure}[!t]
  \centering
  \includegraphics[scale=0.9]{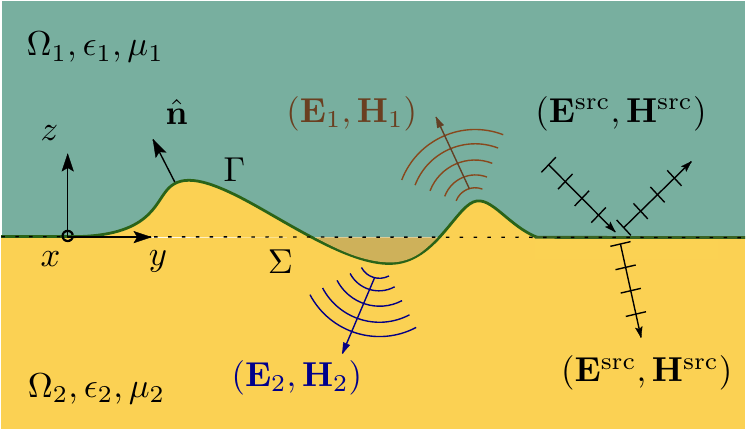}
  \caption{Illustration of the scattering of a plane electromagnetic wave by a locally perturbed penetrable half-space. The boundary $\Gamma$ and the planar surface $\Sigma$ coincide for $(x, y, 0)$ far enough from the bounded local perturbation.}
  \label{fig:scattering}
\end{figure}

In classical layered media SIE formulations~\cite{wannamaker1984electromagnetic,chen2018accurate,yang2011three,Michalski:1997tz,michalski1990electromagnetic,michalski1990electromagneticII}, however, all these attractive features come at the price of employing the dyadic  Green function for layered media, also known as the layer Green function (LGF), which naturally enforces the exact transmission conditions at planar unbounded physical boundaries. The use of the LGF poses difficulties in view of the LGF evaluation cost (there is vast literature on this subject, for which we refer the reader to the  review articles~\cite{Michalski:2016dz,Michalski1998,Aksun:2009fn}). Moreover, problems involving localized perturbations (e.g., open cavities and bumps) give rise to additional challenges to LGF-SIE formulations, as artificial/transparent interfaces need to be introduced in order to properly represent the fields surrounding the perturbations~\cite{perezbruno2014}. 

This work presents a  fully 3D EM layered-media windowed Green function (WGF) method. The WGF method, which was originally developed for scalar layered media problems~\cite{Bruno2015windowed,bruno2017windowed,perez2017windowed} and later extended to waveguides in the frequency and time domains~\cite{garza2020boundary,bruno2017waveguides,sideris2019ultrafast,labarca2019CQ},  completely bypasses the use of  Sommerfeld integrals or other problem-specific Green functions. This is here achieved by first deriving an indirect M\"uller SIE~\cite{muller2013foundations,chao1995regularized} given in terms of free-space Green functions and featuring only weakly-singular kernels, which is posed on the entire unbounded penetrable interface (see Secs.~\ref{sec:inc} and~\ref{sec:SIEFormulation}). The unbounded SIE domain is then effectively truncated to a bounded surface containing the localized perturbations by introducing (in the surface integrals) a smooth windowing function that effectively acts like a reflectionless absorber for the  surface currents leaving the windowed region (see Sec.~\ref{sec:window}). As in the case of the Helmholtz SIEs~\cite{perez2017windowed}, the field errors introduced by the windowing approximation decay faster than any negative power of the diameter of the truncated region. A straightforward (Galerkin) MoM discretization using Rao-Wilton-Glisson (RWG) functions is used to discretize the resulting windowed SIE (see Sec.~\ref{sec:MoM}), although any other Maxwell SIE method could be employed.  A limitation of the proposed approach is that transmission conditions at unbounded penetrable interfaces need to be enforced via second-kind SIEs such as M\"uller's. First-kind SIEs, such as the more popular Poggio--Miller--Chang--Harrington--Wu--Tsai (PMCHWT)~\cite{wu1977scattering,chang1977surface,poggio1970integral}, could be considered  provided they are converted into equivalent second-kind SIEs by means of Calderón preconditioners~\cite{Cools2009}. 

Compared to LGF-SIE formulations, the WGF formulation involves additional unknown surface currents on planar portions of the unbounded dielectric interfaces that eventually lead to larger linear systems. In many cases this additional cost is compensated by the fact that the associated matrix coefficients involve evaluations of the inexpensive free-space Green functions and that the resulting linear system can be efficiently solved iteratively by means of GMRES (see Sec.~\ref{sec:DCIM_comparison}). For problems involving multiple dielectric layers and/or small size PEC inclusions, however, LGF formulations that leverage the discrete complex images method (DCIM)~\cite{yuan2006direct,Alparslan:2010va,ling2000discrete} for the evaluation of the LGF, may well outperform the WGF methodology. 

The proposed approach amounts to a flexible and easy-to-implement MoM for layered media EM problems, in the sense that only minor modifications to existing electromagnetic SIE solvers are needed to deliver the WGF capabilities. The method is thoroughly validated (see Sec.~\ref{sec:examples}) against the exact Mie series scattering solution for a hemispherical bump on a perfectly electrically conducting (PEC) half-space (using a windowed MFIE formulation), and also against the open-source LGF code~\cite{panasyuk2009new}. A performance comparison against a LGF-MoM based on the state-of-the-art C++ library~\emph{Strata}~\cite{sharma2021strata} is presented in Sec.~\ref{sec:DCIM_comparison}. Finally, the proposed methodology is showcased by means of a variety of challenging examples including EM scattering by a large open cavity~(Sec.~\ref{sec:cavity}), an all-dielectric metasurface~(Sec.~\ref{sec:metalens}), and a three-layer plasmonic solar cell  (Sec.~\ref{sec:plasmon}).

\section{Layered media scattering}\label{sec:layered_media}
We consider here the problem of time-harmonic electromagnetic scattering of an incident field $(\elf^\inc,\mgf^\inc)$ by a penetrable locally perturbed half-space~$\Omega_2$, with boundary $\Gamma=\partial\Omega_2$, as depicted in Fig.~\ref{fig:scattering}. Letting $\Omega_1=\R^3\setminus\overline{\Omega_2}$, we express the total electromagnetic field~as
\begin{equation}\label{eq:fld_splitting}
(\elf^{\rm tot},\mgf^{\rm tot}) = (\elf^{\rm src},\mgf^{\rm src })+
(\elf_j,\mgf_j)\ \ \mbox{in}\ \ \Omega_j
\end{equation}
for $j=1,2.$ The known auxiliary source field $(\elf^{\rm src},\mgf^{\rm src })$ which is given in terms of  the incident field $(\elf^\inc,\mgf^\inc)$ under consideration, is constructed so that the fields $(\elf_j,\mgf_j)$, $j=1,2,$ satisfy the homogeneous Maxwell equations 
\begin{equation}\label{eq:PDE}
\nabla\times  \elf_j  -i\omega\mu_j \mgf_j =\bold 0\  \mbox{and}\  \nabla\times  \mgf_j + i\omega \epsilon_j \elf_j =\bold 0\ \mbox{in}\ \Omega_j
\end{equation}
for $j=1,2$, where $\omega>0$ is the angular frequency,  and $\epsilon_j$ and $\mu_j$ are  respectively the permittivity and the permeability within the subdomain~$\Omega_j$. (We have assumed here that the time dependence of the EM fields is given by $\e^{-i\omega t}$.)  For planewave incidences, for instance,  $(\elf^{\rm src},\mgf^{\rm src})$ is taken as the exact total field solution of the problem of scattering of the planewave by the flat lower half-space with planar boundary $\Sigma=\{(x,y,0)\in\R^3\}$ and constants $\epsilon_2$ and $\mu_2$ (see Fig.~\ref{fig:scattering}). The rationale for introducing $(\elf^{\rm src},\mgf^{\rm src})$ lies in ensuring that $(\elf_j,\mgf_j)$, $j=1,2,$ are outgoing wavefields propagating away from the localized perturbations or, more precisely, that they satisfy the Silver-M\"uller radiation condition: 
\begin{equation}
\lim_{|\ner|\to\infty}\left(\sqrt{\mu_j}\,\mgf_j\times\ner-|\ner|\sqrt{\epsilon_j}\,\elf_j\right)=\bold0\ \ \text{in}\ \ \Omega_j,\  j=1,2,
 \label{eq:scattered_maxwell_3D_RC}
\end{equation}
uniformly in all directions $\ner/|\ner|$. The explicit expressions of the source fields utilized throughout the paper are provided in Sec.~\ref{sec:inc} below.

The transmission conditions at the material interfaces, meaning that the tangential components of $(\elf^{\rm tot},\mgf^{\rm tot})$ are continuous across $\Gamma$, lead to the jump conditions
\begin{subequations}\begin{align}
\hat{\mathbf{n}}\times\left\{\mathbf{E}_2|_{-}-\mathbf{E}_1|_{+}\right\} =&~\bol M^{\rm src}\\
\hat{\mathbf{n}} \times\left\{\mathbf{H}_2|_{-}-\mathbf{H}_1|_{+}\right\}=&~\bol J^{\rm src}
\end{align}\label{eq:TC}\end{subequations}
on $\Gamma$, with
\begin{subequations}\begin{align}\label{eq:currentSrc}
\bol M^{\rm src} :=&~\hat{\mathbf{n}} \times\left\{\mathbf{E}^{\rm src}|_{+}-\mathbf{E}^{\rm src}|_-\right\}\\
\bol J^{\rm src}  :=&~\hat{\mathbf{n}} \times\left\{\mathbf{H}^{\rm src}|_{+}-\mathbf{H}^{\rm src}|_-\right\}
\end{align}\label{eq:current_sources}\end{subequations} where we have adopted the notation $\bold F(\ner)|_{\pm} =\lim_{\delta \to 0+}\bold F(\ner\pm\delta\nor(\ner))$  for $\ner\in\Gamma$. As usual, the unit normal vector at $\ner\in\Gamma$ is denoted as~$\nor(\ner)$ and is assumed directed from $\Omega_2$ to $\Omega_1$ (see Fig.1). Existence and uniqueness of solutions of the resulting EM transmission problem are established in~\cite{Cutzach:1998dr}.


\section{Incident and source fields}\label{sec:inc}  Two types of incident fields $(\elf^\inc,\mgf^\inc)$ and corresponding auxiliary source fields $(\elf^{\rm src},\mgf^{\rm src})$ are considered in this paper, namely planewaves and electric dipoles.

Upon impinging on the planar surface $\Sigma$ at the interface between the half spaces $D_1=\{z>0\}$ and $D_2=\{z<0\}$ with wavenumbers $k_1$ and $k_2$ {($k_j=\omega\sqrt{\mu_j\epsilon_j}$, for $j=1, 2$)}, respectively, the incident planewave 
\begin{equation}
\mathbf{E}^{\mathrm{inc}}(\ner)=(\mathbf{p} \times \mathbf{k}) \mathrm{e}^{i \mathbf{k} \cdot \bold{r}} \ \ \text {and} \ \ \mathbf{H}^{\mathrm{inc}}(\ner)=\frac{1}{\omega\mu_1} \mathbf{k} \times \mathbf{E}^{\mathrm{inc}}(\ner)
\label{eq:PW}
\end{equation}
with $\bold p =(p_x,p_y,p_z)$ and $\mathbf{k}=\left(0,k_{1 y},-k_{1 z}\right)$ where $k_{1 z} \geq 0$ and  $|\bold{k}|=\sqrt{k_{1 y}^{2}+k_{1 z}^{2}}=k_{1}$, gives rise to a reflected field $(\elf^{\rm ref},\mgf^{\rm ref})$ in $D_1$ and a transmitted field  $(\elf^{\rm trs},\mgf^{\rm trs})$ in $D_2$. The resulting $x$-independent total field, given by $(\elf^\inc+\elf^{\rm ref},\mgf^\inc+\mgf^{\rm ref})$ in $D_1$ and $(\elf^{\rm trs},\mgf^{\rm trs})$ in $D_2$, is completely determined by the transverse component of the fields~\cite{Chew1995waves}: 
\begin{align*}
\bracenom{E^{\rm inc}_{x}(\ner)}{H^{\rm inc}_{x}(\ner)}=&~\bracenom{E_{0}}{H_{0}}\exp(ik_{1 y} y-ik_{1z}z)\\
\bracenom{E^{\rm ref}_{x}(\ner)}{H^{\rm ref}_{x}(\ner)}=&~\bracenom{E_{0} R^{\mathrm{TE}}}{H_{0} R^{\mathrm{TM}}}\exp(ik_{1 y} y+ik_{1z}z)\\
\bracenom{E^{\rm trs}_{x}(\ner)}{H^{\rm trs}_{x}(\ner)}=&~\bracenom{E_{0} T^{\mathrm{TE}}}{H_{0} T^{\mathrm{TM}}}\exp(ik_{2 y} y-ik_{2z}z)
\end{align*}
depending on the  reflection coefficients:
\begin{align*}
R^{\mathrm{TE}}=\frac{\mu_{2}k_{1z}-\mu_1k_{2z}}{\mu_2k_{1z}+\mu_1k_{2z}},\quad
 R^{\mathrm{TM}}=\frac{\epsilon_2k_{1z}-\epsilon_1k_{2z}}{\epsilon_2k_{1z}+\epsilon_1k_{2z}} 
\end{align*}
the transmission coefficients:
\begin{align*}
 T^{\mathrm{TE}}=\frac{2 \mu_2k_{1z}}{\mu_2k_{1z}+ \mu_1k_{2z}}, \quad T^{\mathrm{TM}}=\frac{2 \epsilon_2k_{1z}}{\epsilon_2k_{1z}+\epsilon_1k_{2z}}
\end{align*}
the amplitudes:  
$$
E_{0}=-p_{z} k_{1 y}-p_{y} k_{1 z},\quad H_{0}=\frac{k_1^2}{\omega\mu_1}p_x
$$ 
and the propagation constants $k_{2y} =k_{1y}$ and $k_{2z} = \sqrt{k_2^2-k^2_{2y}}$ with the complex square root defined so that $\imag k_{2z}\geq 0$. The EM field can be retrieved from the transverse components via
\begin{align*}
\mathbf{E}=&~E_{x} \hat\nex-\frac{1}{i\omega\epsilon} \frac{\partial H_{x}}{\partial z} \hat\ney+\frac{1}{i\omega\epsilon} \frac{\partial H_{x}}{\partial y} \hat\nez \\
\mathbf{H}=&~H_{x} \hat\nex+\frac{1}{i\omega\mu} \frac{\partial E_{x}}{\partial z} \hat\ney-\frac{1}{i\omega\mu} \frac{\partial E_{x}}{\partial y} \hat\nez.
\end{align*}
With these expressions at hand we define the auxiliary planewave source field as 
\begin{align}
(\elf^{\rm src},\mgf^{\rm src}) =\begin{cases}(\elf^{\rm inc},\mgf^{\rm inc})+(\elf^{\rm ref},\mgf^{\rm ref})&\text{in }\Omega_1\\
(\elf^{\rm trs},\mgf^{\rm trs})&\text{in }\Omega_2.
\end{cases}\label{eq:PW_src}
\end{align}
Since by  construction the source field~\eqref{eq:PW_src} satisfies the exact transmission conditions at $\Sigma$, it holds that the current sources $\bol M^{\rm src}$ and $\bol J^{\rm src}$ defined in~\eqref{eq:current_sources} are supported on the (bounded) local perturbation~$\Gamma\setminus\Sigma$. 

In the special case when $\Omega_2$ is occupied by a PEC, in which the boundary condition $\nor\times\elf^{\rm tot} =\bol 0$ holds on the interface $\Gamma$, we have that $(\elf^{\rm src},\mgf^{\rm src})$ takes the form~\eqref{eq:PW_src} with $(\elf^{\rm trs},\mgf^{\rm trs})= (\bol 0,\bol 0)$ and $(\elf^{\rm ref},\mgf^{\rm ref})$ given in terms of the reflection coefficients $R^{\TE} =-R^{\TM} = -1$. 

Finally, we take
\begin{equation}
\mgf^\inc = \frac{1}{i\omega\mu_j}\nabla\times\{G_j(\cdot,\ner_0){\bf p}\},\ \ \elf^\inc= \frac{-1}{i\omega \epsilon_j}\nabla\times \mgf^{\inc}\label{eq:dipoleInc}
\end{equation}
with 
\begin{equation}\label{eq:freespaceGF}
G_j(\ner,\ner') := \frac{\e^{ik_j|\ner-\ner'|}}{4\pi|\ner-\ner'|}
\end{equation}being the (Helmholtz) free-space Green function, as the incident field produced by an electric dipole at $\ner_0\in\Omega_j$. The corresponding source field is thus selected as
\begin{align}
(\elf^{\rm src},\mgf^{\rm src}) =\begin{cases}(\elf^{\rm inc},\mgf^{\rm inc})&\text{ in }\Omega_j\\
(\bol 0,\bol 0)&\text{ otherwise}\end{cases}
\label{eq:PT_src}\end{align}
for $j=1,2$. 

\section{Second-kind integral equation formulation\label{sec:SIEFormulation}}
In order to approximate the unknown EM  fields $(\elf_j,\mgf_j)$, $j=1,2$, we resort to a second-kind indirect M\"uller formulation. We start by introducing the off-surface integral operators
\begin{align}
\begin{split}
(\mathcal S_j\bol\varphi)(\ner):=&~\int_{\Gamma}G_j(\ner,\ner')\bol\varphi(\ner') \de s'+\\
&
\frac{1}{k^2_j}\nabla\int_{\Gamma}G_j(\ner,\ner')\ \nabla_{s}'\cdot\bol \varphi (\ner') \de s' \label{eq:EFIE_pot}
\end{split}\\
(\mathcal{D}_{j}\bol\varphi)(\mathbf{r}):=&~\nabla \times \int_{\Gamma} G_{j}(\mathbf{r},\mathbf{r}^{\prime}) \bol\varphi(\mathbf{r}^{\prime})\de s' \label{eq:MFIE_pot}
\end{align}
for $\ner\in\R^3\setminus\Gamma$, with $\bol\varphi$ being a vector field tangential to $\Gamma$.
(In what follows the surface integrals over $\Gamma$ must be interpreted as conditionally convergent.) 
The unknown EM fields $(\elf_j,\mgf_j)$ are thus sought as
\begin{subequations}\begin{align}
\elf_j(\ner) :=&~k_j^2(\mathcal{S}_{j}\bol v)(\ner)+i\omega\mu_j(\mathcal{D}_{j}\bol u)(\ner)\\
\mgf_j(\ner) :=&~k_j^2(\mathcal{S}_{j}\bol u)(\ner)-i\omega\epsilon_j(\mathcal{D}_{j}\bol v)(\ner)
\end{align}\label{eq:SC}\end{subequations}
for $\ner\in\Omega_j$, $j=1,2$, in terms of unknown surface currents~$\bol u$ and $\bol v$ that are to be determined by means of a SIE posed on  $\Gamma$.  Clearly, the field defined in~\eqref{eq:SC} satisfy Maxwell equations~\eqref{eq:PDE}, in view of the fact that $\nabla \times \mathcal{S}_j = \mathcal{D}_j$ and $\nabla \times \mathcal{D}_j = k_j^2 \mathcal{S}_j$.

In order to derive a SIE for the currents we make use of the well-known jump relations:
\begin{align}
  \nor\times (\mathcal S_j\bol \varphi)|_{\pm} =\mathcal  T_j\bol \varphi \ \text{ and }\ 
 \nor\times (\mathcal D_j\bol \varphi)|_{\pm} =\mathcal  K_j\bol \varphi\pm\frac{\bol\varphi}{2}
\label{eq:jumps} \end{align} on $\ner\in\Gamma$,  where 
\begin{align}
  (\mathcal T_j\bol \varphi)(\ner):=\,&\nor(\ner)\times\int_{\Gamma}G_j(\ner,\ner')\bol\varphi(\ner') \de s'+\nonumber\\
  &\frac{1}{k_j^{2}}\nor(\ner)\times\nabla\int_{\Gamma}G_j(\ner,\ner')\ \nabla_{s}'\cdot\bol \varphi (\ner') \de s'\label{eq:E2}
\end{align}  and
\begin{align}\label{eq:Kop}
(\mathcal K_j\bol\varphi)(\ner):=\,&\nor(\ner)\times\nabla\times\int_{\Gamma}G_j(\ner,\ner')\bol \varphi(\ner') \de s'.
\end{align}  

Evaluating the integral representation formulae~\eqref{eq:SC} on $\Gamma$ and using~\eqref{eq:jumps} we obtain
\begin{subequations}\begin{align}
\nor\times\elf_j|_\pm =k_j^2\mathcal{T}_{j}\bol v+i\omega\mu_j\lf\{\pm\frac{\bol u}{2}+\mathcal{K}_{j}\bol u\rg\}
\end{align}
for the electric fields, and
\begin{align}
\nor\times\mgf_j|_\pm =k_j^2\mathcal{T}_{j}\bol u-i\omega\epsilon_j\lf\{\pm\frac{\bol v}{2}+\mathcal{K}_{j}\bol v\rg\}
\end{align}\label{eq:on_surface}\end{subequations}
for the magnetic fields. Therefore, enforcing the transmission conditions~\eqref{eq:TC} by taking the appropriate linear combination of the relations~\eqref{eq:on_surface}, we arrive at the following  SIE  for the unknown vector of current densities $\lf[\bol u,\bol v\rg]^\top$:
\begin{equation}
\frac{i\omega}{2}\left[\begin{array}{c}
\!\!-(\mu_1+\mu_2)\bol u\!\! \\
(\epsilon_1+\epsilon_2)\bol v\!\!
\end{array}\right]+\mathsf{T}\left[\begin{array}{c}
\!\!\bol u \!\!\\
\!\!\bol v\!\!
\end{array}\right]=\left[\begin{array}{c}
\!\!\bol M^{\rm src} \!\!\\
\!\!\bol J^{\rm src}\!\!
\end{array}\right]\quad\text{on}\quad\Gamma
\label{eq:LMIE}\end{equation}
where the block operator $\mathsf{T}$ is given by 
\begin{equation}
\mathsf{T}=\left[\begin{array}{cc}
i\omega\left(\mu_2\mathcal K_{2}-\mu_1\mathcal K_{1}\right) & k^2_2\mathcal T_{2}-k^2_1\mathcal T_{1} \\
k^2_2\mathcal T_{2}-k^2_1\mathcal T_{1} & -i\omega\left(\epsilon_2\mathcal K_{2}-\epsilon_1\mathcal K_{1}\right)
\end{array}\right].\label{eq:transOp}\end{equation}

We emphasize that the rationale underlying expressing the EM fields as in~\eqref{eq:SC} lies in making the strongly singular operators $\mathcal T_j$, $j=1,2$, appear in the resulting SIE~\eqref{eq:LMIE} as the linear combination $k^2_2\mathcal T_{2}-k^2_1\mathcal T_{1}$. Indeed, as shown in~\cite{yla2005well,chao1995regularized} and in Sec.~\ref{sec:MoM} below, this linear combination can be cast into a bounded integral kernel tractable by standard off-the-shelf quadrature rules.

Finally, it is worth mentioning that the two-layer media scattering problem considered in this section can as well be recast as the classical direct M\"uller integral equation involving the same operator $\mathsf T$, but with a different right-hand-side that entails evaluation of the singular $\mathcal T_1$ operator~\cite[Sec.~6.2]{perez2017windowed}.

\section{Windowed Green function method\label{sec:window}}
The fact that the SIE~\eqref{eq:LMIE} is posed on an unbounded surface~$\Gamma$ introduces the salient issue of having to suitably truncate the computational domain to numerically approximate the SIE solution via the MoM. Therefore, instead of solving~\eqref{eq:LMIE} on the entire material interface~$\Gamma$, we make use of a windowed SIE to obtain approximations of the surface current densities $[\bol u,\bol v]^{\top}$ over the relevant portion of~$\Gamma$ containing the localized perturbations.  In order to do so we introduce a slow-rise infinitely smooth window function $w_{A}:\R^3\to\R$ which vanishes with all its derivatives outside the cylinder $\{\sqrt{x^2+y^2}<A\}\times\R$. More precisely, the window function is selected as $w_A(\ner) = \eta(\sqrt{x^2+y^2},cA,A)$  for $A>0$, $0<c<1$, and
\begin{equation}\begin{split}&\eta\left(s, s_{0}, s_{1}\right):=\\ &\left\{\begin{array}{cl}
1 & \text { if }|s|<s_{0} \\
\displaystyle\exp \left(\frac{2 \mathrm{e}^{-1 / b}}{b-1}\right), b=\frac{|s|-s_{0}}{s_{1}-s_{0}} & \text { if } s_{0}<|s|<s_{1} \\
0 & \text { if }|s| \geq s_{1}.
\end{array}\right.\end{split}\end{equation}
The parameter value $c = 0.7$ is used in all the numerical examples presented in Sec.~\ref{sec:examples}. (Other definitions of the window function, such as $w_A(\ner) = \eta(x,cA,A)\eta(y,cA,A)$ or  $w_A(\ner) = \eta(x,cA_x,A_x)\eta(y,cA_y,A_y)$ with $A_x,A_y>0$,  for instance, can also be employed so as to suitably adjust $\widetilde\Gamma_A$ to the particular shape of the localized perturbations.)

We then consider the following windowed~SIE:
\begin{equation}
\frac{i\omega}{2}\left[\begin{array}{c}
\!\!-(\mu_1+\mu_2){\bol u}_{\!A}\!\! \\
(\epsilon_1+\epsilon_2){\bol v}_{\!A}\!\!
\end{array}\right]+\mathsf{T}_{\!A}\left[\begin{array}{c}
\!\!{\bol u}_{\!A} \!\!\\
\!\!{\bol v}_{\!A}\!\!
\end{array}\right]=\left[\begin{array}{c}
\!\!\bol M^{\rm src}\!\!\\
\!\!\bol J ^{\rm src}\!\!
\end{array}\right]\ \ \text{on}\ \ \Gamma_A\label{eq:WLMIE}
\end{equation}
 where $ \Gamma_{\!A}=\left\{\ner \in \Gamma: w_{\!A}(\ner) \neq 0\right\}$ and where the windowed operator $\mathsf{T}_{\!A}$ is defined as $\mathsf{T}$ in~\eqref{eq:transOp} but with the \emph{windowed Green function}
\begin{equation}G_{A,j}(\ner,\ner') = w_A(\ner')G_j(\ner,\ner')\label{eq:winGF}\end{equation}
replacing the free-space Green function $G_j$ appearing in the definition of $\mathcal T_j$ and $\mathcal K_j$ in~\eqref{eq:E2} and~\eqref{eq:Kop}, respectively. 

Existence and uniqueness of  solutions of the windowed SIE~\eqref{eq:WLMIE} can be established (under reasonable smoothness assumptions on $\Gamma_A$ and up to a  countable set of frequencies~$\omega$) by invoking the Fredholm alternative, which holds true in this case by virtue of the compactness of $\mathsf{T}_{\!A}$ (in an appropriate function space). Alternatively, for sufficiently small contrasts $\epsilon_1/\epsilon_2$ and $\mu_1/\mu_2$, existence and uniqueness could be established following a Neumann series approach.
 
 As it turns out, $[\bol u_{\!A},\bol v_{\!A}]^\top$ provides an excellent approximation of the exact currents $[\bol u,\bol v]^\top$ within $\widetilde\Gamma_{\!A} =\{\ner\in\Gamma: w_{\!A}(\ner)=1\}$. Indeed, as in the two-dimensional electromagnetic case~\cite{Bruno2015windowed,perez2017windowed,bruno2017waveguides,sideris2019ultrafast}, we have that the errors in the approximation $[\bol u,\bol v]^\top \approx [\bol u_{\!A},\bol v_{\!A}]^\top$ decay super-algebraically fast in $\widetilde\Gamma_{A_0}$ for a fixed $A_0>0$, as the window size $A$ increases.

With the surface current densities $[\bol u_{A},\bol v_{A}]^{T}$ at hand, the approximate EM  fields can be easily obtained by, respectively, substituting $\bol u$ and $\bol v$ by $\bol u_{A}$ and $\bol v_{A}$ in the representation formula~\eqref{eq:SC}, and by replacing $G_j$ by the WGF~\eqref{eq:winGF} in the off-surface operators $\mathcal S_j$ and $\mathcal D_j$ defined in~\eqref{eq:EFIE_pot} and~\eqref{eq:MFIE_pot}, respectively. These substitutions produce the approximate fields
\begin{subequations}\begin{align}
\widetilde\elf_{j}(\ner) :=&k_j^2(\mathcal{S}_{A,j}{\bol v}_{\!A})(\ner)+i\omega\mu_j(\mathcal{D}_{A,j}{\bol u}_A)(\ner)\label{eq:windElf}\\
\widetilde\mgf_{j}(\ner) :=&k_j^2(\mathcal{S}_{A,j}{\bol u}_{\!A})(\ner)-i\omega\epsilon_j(\mathcal{D}_{A,j}{\bol v}_{\!A})(\ner)
\end{align}\label{eq:windEMfld}\end{subequations}
for $\ner\in\Omega_j$, $j=1,2$, where $\mathcal S_{A,j}$ and $\mathcal D_{A,j}$ are the resulting windowed off-surface operators.

We note that although formula~\eqref{eq:windEMfld} does not directly yield accurate far-fields, they can still be obtained from the accurate near-fields produced by~\eqref{eq:windEMfld} within $\{\ner\in\R^3: w_A(\ner)=1\}$ in a manner akin to~\cite[Sec.~3.6]{Bruno2015windowed} for the corresponding scalar problem (which in this case would involve the leading-order asymptotic approximation of the dyadic  LGF, $\bold G(\ner,\ner')$, as $|\ner|\to\infty$~\cite[Sec.~2.6]{Chew1995waves}).

As mentioned above, a limitation of the proposed approach is that first-kind SIEs, such as the PMCHWT formulation~\cite{wu1977scattering,chang1977surface,poggio1970integral}, which  has been the preferred approach for electromagnetic transmission problems~\cite{zhu2004comparison}, is not directly compatible with the WGF approach. In a nutshell, windowed kernels decay exponential fast near the edges of the truncated surface, thus allowing surface currents near those edges to lie in the approximate nullspace of the PMCHWT WGF-MoM matrices. Such matrices have then eigenvalues very close to the origin making the linear system too ill-conditioned to be accurately solved by either direct or iterative methods. In contrast, second-kind SIEs, like the ones employed in this contribution, do not suffer from this problem because the identity term shifts the spectrum sufficiently far away from the origin (see Fig. \ref{fig:eigen} in Sec.~\ref{sec:SomProbblem}). A possible remedy to this issue is the use of Calderón preconditioners~\cite{Cools2009}, which take advantage of the operators' self-regularizing properties to convert  first-kind SIEs into equivalent well-conditioned second-kind SIEs.

\section{MoM discretization}\label{sec:MoM}
We start off this section by considering a triangulation of the truncated surface $\Gamma_{\!A}$ which lies within the support of the window function $w_{\!A}$.
 We then expand $[\bol u_A,\bol v_A]^\top$  in terms of the standard div-conforming RWG basis functions~\cite{rao1982electromagnetic}. In detail, we let
\begin{align}
\bol u_A(\ner) \approx \sum_{n=1}^N u_n\bold f_n(\ner)\andtext \bol v_A(\ner) \approx \sum_{n=1}^N v_n\bold f_n(\ner)\label{eq:current_exp}
\end{align}
for $\ner\in\Gamma_A$, where $N$ is the total number of edges in the triangular mesh, and $\bold f_n$ are the RWG basis functions~\cite{gibson2014method}.




As in~\cite{yla2005well}, we apply the Galerkin scheme to determine the coefficients $u_n$ and $v_n$ in the approximations~\eqref{eq:current_exp} by replacing~\eqref{eq:current_exp} in the windowed SIE~\eqref{eq:WLMIE} and then testing the resulting equations against the same div-conforming basis functions $\bold f_n$. 
We thus obtain the following linear system for the coefficients:
\begin{equation}\label{eq:MoMLS}
M{\bf x} = {\bf b} 
\end{equation}
where ${\bf x} = [u_1,\ldots,u_N,v_1,\ldots,v_N]^\top\in\mathbb C^{2N}$,
\begin{equation}
M = \lf[\begin{array}{ccc}
M^{(1,1)} & M^{(1,2)}\\
M^{(2,1)} & M^{(2,2)}
\end{array}\rg]\in\mathbb C^{2N\times 2N}
\label{eq:sys_matrix}\end{equation}
with blocks $M^{(p,q)}$, $p,q=1,2$ defined as 
$$
M^{(p,q)}_{m,n} = \int_{\Gamma_A}\bold f_m(\ner)\cdot (\mathcal M^{(p,q)}\bold f_n)(\ner)\de s,\quad n,m=1\ldots,N,
$$
in terms of the operators:
\begin{equation*}\begin{split}
\mathcal M^{(1,1)}\bol\varphi  =\!& 
-\frac{i\omega(\mu_1+\mu_2)}{2}\bol\varphi +i\omega\left(\mu_2\mathcal K_{A, 2}-\mu_1\mathcal K_{A, 1}\right)[\bol\varphi]\smallskip\\
\mathcal M^{(1,2)}\bol\varphi=&\,\mathcal M^{(2,1)}\bol\varphi =(k^2_2\mathcal T_{A, 2}-k^2_1\mathcal T_{A, 1})[\bol\varphi]\smallskip\\
\mathcal M^{(2,2)}\bol\varphi =& \frac{i\omega(\epsilon_1+\epsilon_2)}{2}\bol\varphi+i\omega\left(\epsilon_1\mathcal K_{A, 1}-\epsilon_2\mathcal K_{A, 2}\right)[\bol\varphi].
\end{split}\end{equation*}
On the other hand, the right-hand-side vector ${\bf b}\in \mathbb C^{2N}$ is given by ${\bf b}_m = \int_{\Gamma_A}{\bf f}_m\cdot {\bf M^{\rm src}}\de s$ for $m=1,\ldots,N$, and ${\bf b}_m = \int_{\Gamma_A}{\bf f}_m\cdot {\bf J^{\rm src}}\de s$ for $m=N+1,\ldots,2N$.

Evaluation of the matrix entries boils down to compute integrals of the form
\begin{equation*}
  I^{(1)}_{m,n}
=\int_{\Gamma_{\!A}}\!\!{\bf f}_m (\ner)\cdot\lf\{\nor(\ner)\times\nabla\times\int_{\Gamma_{\!A}}\!\!G_{\!A,j}(\ner,\ner'){\bf f}_n (\ner') \de s'\rg\}\de s
\end{equation*}
in the case of diagonal blocks $M^{(1,1)}$ and $M^{(2,2)}$, and integrals of the form
\begin{eqnarray*}
\hspace{-1cm}&&I^{(2)}_{m,n}
=\int_{\Gamma_{\!A}}{\bf f}_m(\ner)\cdot 
\lf\{\nor(\ner)\times \int_{\Gamma_{\!A}}G_{\!A,j}(\ner,\ner')\bold f_n (\ner') \de s'\rg\}\de s\label{eq:Z2}\\
\hspace{-1cm}&&\begin{split}I^{(3)}_{m,n}=&\int_{\Gamma_{\!A}}{\bf f}_m(\ner)\cdot 
\bigg\{\nor(\ner)\;\times\\& \int_{\Gamma_{\!A}}\nabla\lf[G_{\!A,2}(\ner,\ner')-G_{\!A,1}(\ner,\ner')\rg]\nabla_s'\cdot\bold f_n (\ner') \de s'\bigg\}\de s,\end{split}\label{eq:Z3}
\end{eqnarray*}
in the case of the off-diagonal blocks $M^{(1,2)}$ and $M^{(2,1)}$, where $G_{\!A,j}$, $j=1,2,$ are defined in~\eqref{eq:winGF}. Note that the integrands above become singular whenever the supports of~${\bf f}_n$ and ${\bf f}_m$ intercept each other. However, given that
$$
\nabla [G_2(\ner,\ner')-G_1(\ner,\ner')] =-\frac{(k_2^2 - k_1^2)}{2}\frac{(\ner-\ner')}{|\ner-\ner'|} + o(1)
$$
as $|\ner-\ner'|\to 0$, we have that the integrand in the definition of $I_{m,n}^{(3)}$ remains bounded. In the numerical examples considered in the next section, we utilize the Duffy-like singularity cancellation technique presented in~\cite{Sauter2010} to render these weakly-singular integrands into smooth functions that we integrate by means of standard Gauss quadrature rules.

Finally, we briefly discuss the selection of the parameter $A$ and the mesh size $h$ associated with the discretization of $\Gamma_A$. Since the WGF truncation errors in $\mathcal T_{A,j}$ and $\mathcal K_{A,j}$, $j=1,2$,  decay faster than any power of $|k_jA|^{-1}$ as $A$ increases~\cite{perez2017windowed}, $A$ should be selected such that $A>\max\{k_1^{-1},|k_2|^{-1}\}$. It was found in practice that  $A>16\pi\max\{k_1^{-1},|k_2|^{-1}\}+R$, where $R>0$ is the radius of the smallest ball containing the perturbations, is more than enough to suppress any error stemming from the windowing approximation, making the overall WGF-MoM error of order $O(h^2)$ for any reasonable small mesh size $h$ (see Figs.~\ref{fig:errorPEC} and~\ref{fig:errorDipole}). Regarding the discretization of $\Gamma_A$, on the order hand, it has to be such that the spatial oscillations of the surface integrands in the integral operators are well resolved, i.e., $h<\pi\max\{ k_1^{-1},(\real{k_2})^{-1}\}$ so that  Nyquist criterion is not violated.

\section{Validation and examples\label{sec:examples}} 
A variety of numerical examples are presented in this section to validate and demonstrate  the accuracy, efficiency and applicability of the proposed methodology.

\subsection{PEC hemispherical bump} \label{sec:bump_example}
First, in order to validate the proposed WGF-MoM approach, we  consider the problem of scattering of an incident EM planewave~\eqref{eq:PW} of  wavelength $\lambda=2\pi/k_1=1$, by a PEC hemispherical bump of radius $\lambda$ placed on top of the PEC half-space $\{z<0\}$ (see inset in Fig.~\ref{fig:errorPEC}).  We thus compare the exact solution $\elf^{\rm ref}$~(see Appendix~\ref{app:A})  with the numerical WGF-MoM solution, which is obtained using an indirect windowed MFIE formulation.  We note that since the exact PEC half-space Green function can be computed in closed form (via the method of images), this as well as more general PEC obstacle and bump-like scattering problems can be directly cast into the classical MFIE and EFIE posed on the obstacle/bump's  surface. There is therefore no particular advantage of employing the WGF method in these cases. The reason why this problem is here considered is that, to the best of the authors' knowledge, this is the only problem of scattering by a locally perturbed infinite planar surface that admits an exact Mie series solution.

In detail, the approximate total electric field takes the form $\widetilde\elf^{\rm tot} = \widetilde\elf^s + \elf^{\rm src}$, where $\elf^{\rm src}$ is given in Sec.~\ref{sec:inc} and $\widetilde\elf^s = \mathcal D_{A, 1}[\bol u_A]$ with the currents $\bol u_A$ being the solution of the windowed MFIE:
$$
\frac{\bol u_{\!A}}{2} + \mathcal K_{A, 1}[\bol u_{\!A}] = -\nor\times\elf^{\rm src}\quad\mbox{on}\quad\Gamma_{\!A}, 
$$
which is solved using the standard MoM discretization~\cite{gibson2014method}. (The corresponding total magnetic field can be retrieved from $\bol u_A$ via $\widetilde\mgf^{\rm tot} = \widetilde\mgf^s + \mgf^{\rm src}$ with $\widetilde\mgf^s = \mathcal -i\omega\epsilon_1 \mathcal S_1[w_A\bol u_A]$.)

Figure~\ref{fig:errorPEC} displays the electric field errors 
\begin{equation}
{\rm error }=\max _{\ner \in \Theta}\left|\widetilde\elf^{\rm tot}(\ner)-\mathbf{E}^{\mathrm{ref}}(\ner)\right| / \max _{\ner \in \Theta}\left|\mathbf{E}^{\mathrm{ref}}(\ner)\right|\label{eq:rel_error}
\end{equation}
obtained by evaluating both solutions at a fixed target point set~$\Theta$ on a hemispherical surface of radius $2\lambda$ concentric to the PEC bump, for various mesh ($h>0$) and window ($A>0$) sizes. TE- and TM-polarized planewave incident fields~\eqref{eq:PW}  at the grazing angle $\pi/32=\arctan(k_{1z}/k_{1y})$ were used in these examples. The respective linear systems~\eqref{eq:MoMLS} were iteratively solved by means of GMRES~\cite{saad1986gmres} which converged in about 20 iterations using a relative tolerance of  $10^{-5}$ and the Jacobi (diagonal) preconditioner. Almost identical results are obtained for incidences closer to normal.

\begin{figure}[!t]
    \centering
    \includegraphics[scale=0.55]{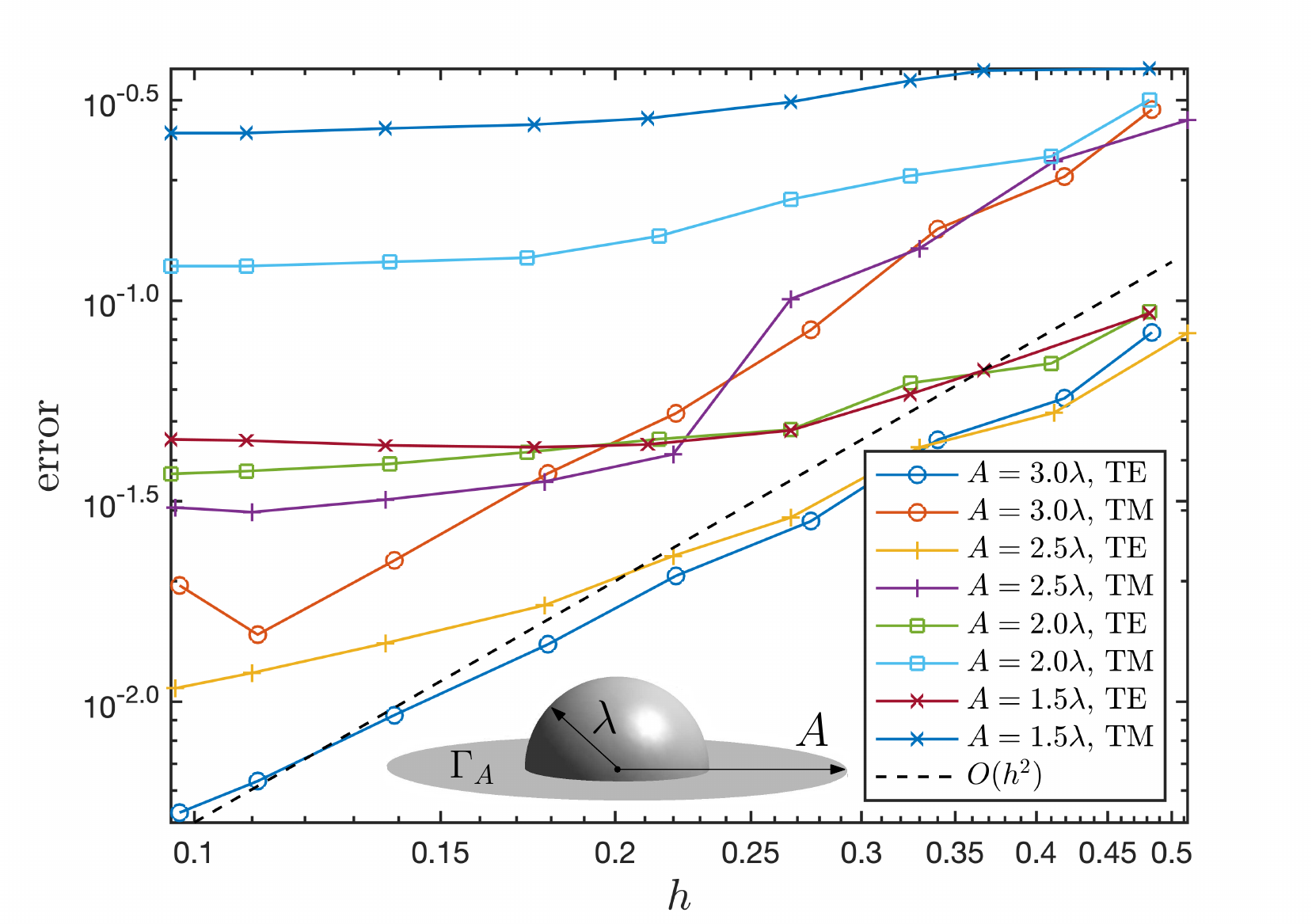}
    \caption{Errors~\eqref{eq:rel_error} in the windowed MFIE-MoM solution of the PEC hemispherical bump problem for various approximately uniform meshes of size $h>0$ and window sizes $A>0$, plotted in log-log scale. The incident fields correspond to a planewave with grazing angle of $\pi/32$ in both TE and TM polarizations. The dashed line marks the second-order slope.}
    \label{fig:errorPEC}
    \end{figure}	

There are two types of errors present in these results; the error stemming from the WGF approximation, which decreases super-algebraically as $A\to\infty$, and the MoM error, which decreases as $h^2$ as $h\to 0$. The former becomes dominant for small $A$ values, as can be seen in the flattening of the error curves for small $h$ values, while the latter becomes dominant for sufficiently large $A$ values, as can be seen in the quadratic error decay as~$h$ decreases. These results validate the convergence of the our windowed MoM-solver, which is not affected by the planewave incidence angle and polarization.

We mention in passing  that this simple windowed MFIE formulation can as well be employed to tackle the rather classical PEC open cavity problem~(cf.~\cite[Ch.~10]{jin2015finite}). Standard SIE formulations for this problem~\cite{ammari2000integral,wood1999development,PerezArancibia:2014fg} entail introducing an artificial transparent surface to close the open cavity, which is not  needed by the windowed MFIE formulation. 

\begin{figure}[!t]
    \centering
    \includegraphics[scale=0.55]{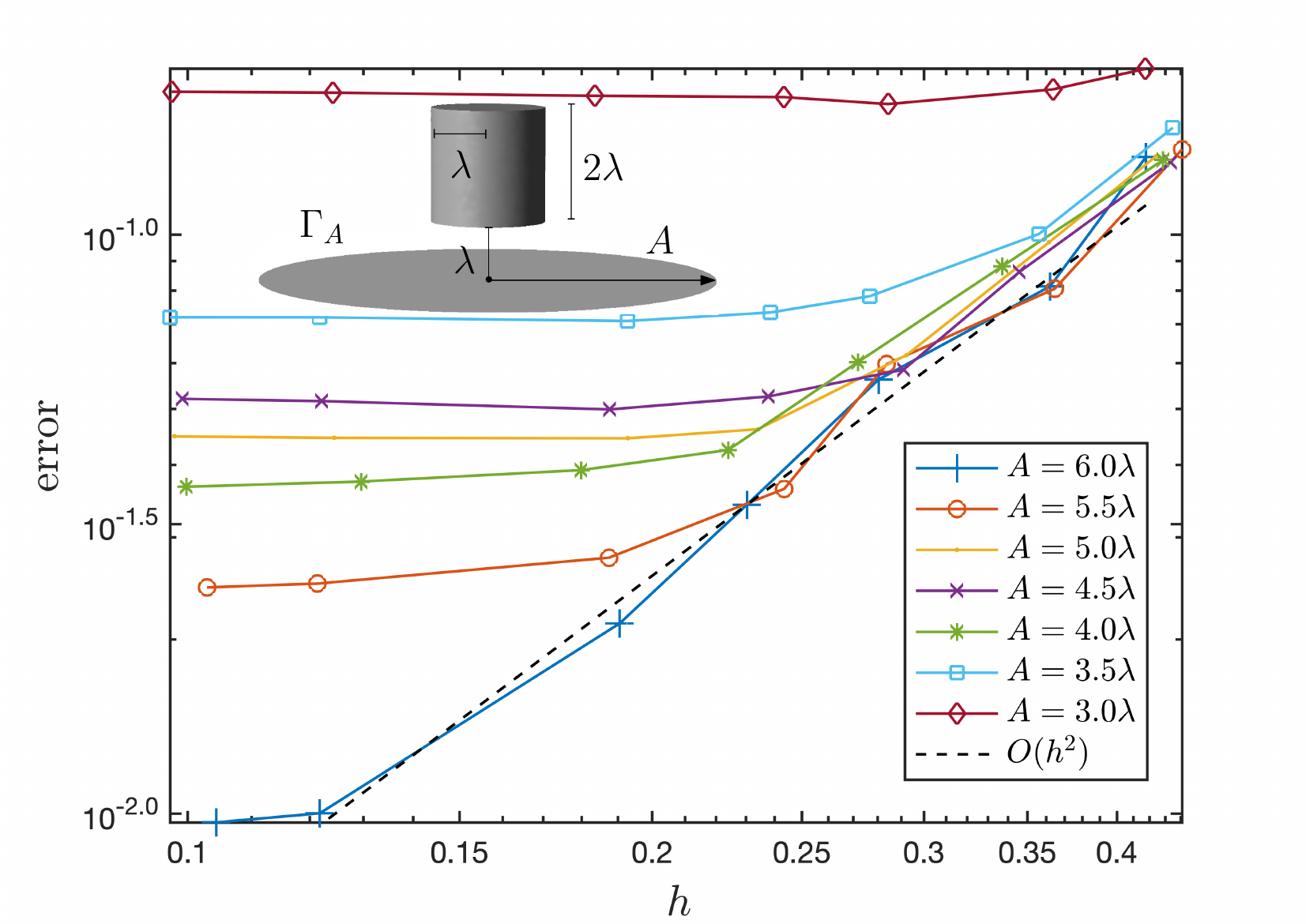}
    \caption{Errors~\eqref{eq:rel_error} in the WGF-MoM solution of the Sommerfeld half-space problem for various approximately uniform meshes of size $h>0$ and window sizes $A>0$, plotted in log-log scale. Ten electric dipole sources uniformly randomly placed within the boundary of the cylinder shown in the inset figure, were used as the incident field. The dashed line marks the second-order slope.}
    \label{fig:errorDipole}
\end{figure}

 \subsection{Sommerfeld half-space problem}\label{sec:SomProbblem}
Our next example deals with the classical Sommerfeld half-space problem~\cite{Michalski:2016kk}. We consider an incident electric field $\elf^\inc$ produced by the superposition of ten randomly placed (at the points $\ner_\ell$, $\ell=1,\ldots,10$) electric dipoles sources~\eqref{eq:dipoleInc}. The dipole sources are uniformly distributed within the boundaries of a cylinder of radius $\lambda$, height $2\lambda$, and centered at $(0,0,2\lambda)$ (see inset in Fig.~\ref{fig:errorDipole}), where $\lambda=2\pi/k_1=1$~m in this case. The incident field  impinges on a dielectric half-space $\Omega_2=\{z<0\}$ with $k_2=\sqrt{2}k_1$ ($\epsilon_2=2\epsilon_1$). The exact total electric field takes the form $\elf^{\rm ref}(\ner) = \sum_{\ell=1}^{10}\bold G(\ner,\ner_\ell)\bold p_\ell$ where~$\bold G$ is the dyadic LGF~\cite{Chew1995waves} and~$\bold p_\ell$, $\ell=1,\ldots,10$, are random polarization unit vectors. 

The approximate total electric field, on the other hand, is given by $\widetilde\elf^{\rm tot} =\widetilde \elf_1+\elf^{\rm src}$ in~$\Omega_1$, where~$\widetilde \elf_1$ is obtained from~\eqref{eq:windElf} with currents ($\bol u_A,\bol v_A$) produced by means of the MoM applied to the windowed SIE~\eqref{eq:WLMIE}, and where the source field $\elf^{\rm src}$ is given in~\eqref{eq:PT_src}.

\begin{table}[h!]
	\centering
	\caption{Errors in Fig. \ref{fig:errorDipole} for the seven window sizes $A$ used, which are indexed by $n=1,\dots 7$, corresponding to $h\approx 0.1$. The log-log slope of the error (as a function of $A$) is computed as $\sigma_n=-\log\left({\rm error}_n / {\rm error}_{n-1}\right) / \log\left( A_n /A_{n-1}\right)$ for $n=2,\dots,7$. The increasing  $\sigma_n$ values for $n\geq 4$ demonstrate the super-algebraic convergence of the WGF-MoM achieved for sufficiently small mesh sizes (when the windowing error is dominant).} 
	\resizebox{\columnwidth}{!}{%
	\begin{tabular}{c|c|c|c|c|c|c|c}
		\toprule	
		$n$             & {1} & {2} & {3} & {4} & {5} & {6} & {7} \\\hline
		$ A_n / \lambda$   & 3          & 3.5        & 4          & 4.5        & 5          & 5.5        & 6       \\   
		${\rm error}_n\times10^{2}$         & $17$     & 7.2     & 5.2     & 4.4     & 3.6     & 2.4     & 1.0     \\ 
		$\sigma_n$  &    $-$       & $5.8$       & $2.4$       & $1.2$       & $1.9$       & $4.2$       & $10.7$      \\ \bottomrule
	\end{tabular}
	}
	\label{tab:win_vs_error}
\end{table}

Figure~\ref{fig:errorDipole} displays the errors~\eqref{eq:rel_error} in the total field~$\widetilde\elf^{\rm tot}$ for various mesh and window sizes. The target point set $\Theta$ used to compute the errors encompasses 2,332  points lying on the surface of the cylinder containing the dipole sources (see inset in Fig.~\ref{fig:errorDipole}). The reference (total) field in this example, was produced by the LGF FORTRAN~code~\cite{panasyuk2009new}. The particular source-target point configuration of Fig.~\ref{fig:errorDipole} intentionally avoids dealing with difficult cases that could affect the accuracy of the LGF evaluations.

Once again, fast convergence is observed as the window size $A>0$ increases while the expected second-order convergence is attained as the mesh size $h>0$ decreases. The smallest errors reported in Fig. \ref{fig:errorDipole} (for $h\approx 0.1$) are reproduced in Table~\ref{tab:win_vs_error} for the various window sizes used in this example. In view of the fact that the log-log slope $\sigma_n$ grows (in magnitude) as $A$ increases, we have that the error decays super-algebraically fast as $A$ increases (algebraic convergence of any fixed order would produce an approximately constant slope $\sigma_n$).

\begin{figure}[!t]
    \centering
    \includegraphics[width=4.35cm]{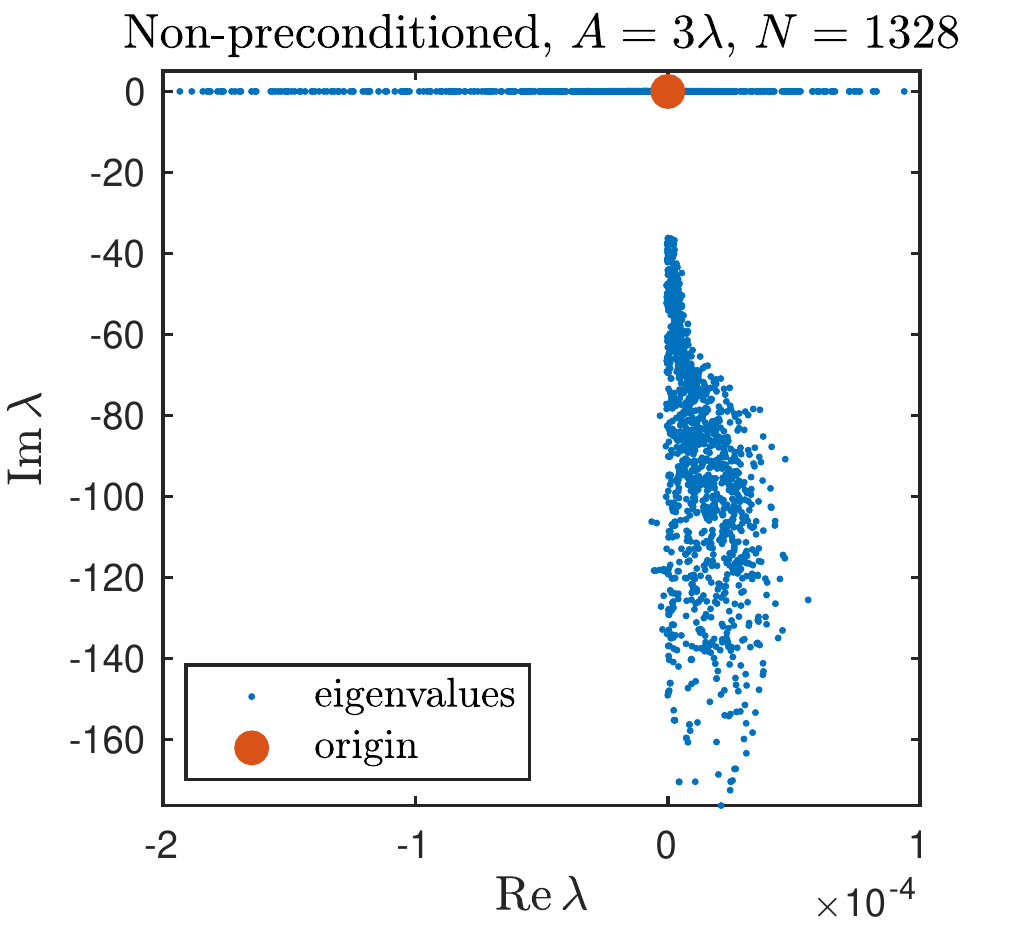}
    \includegraphics[width=4.35cm]{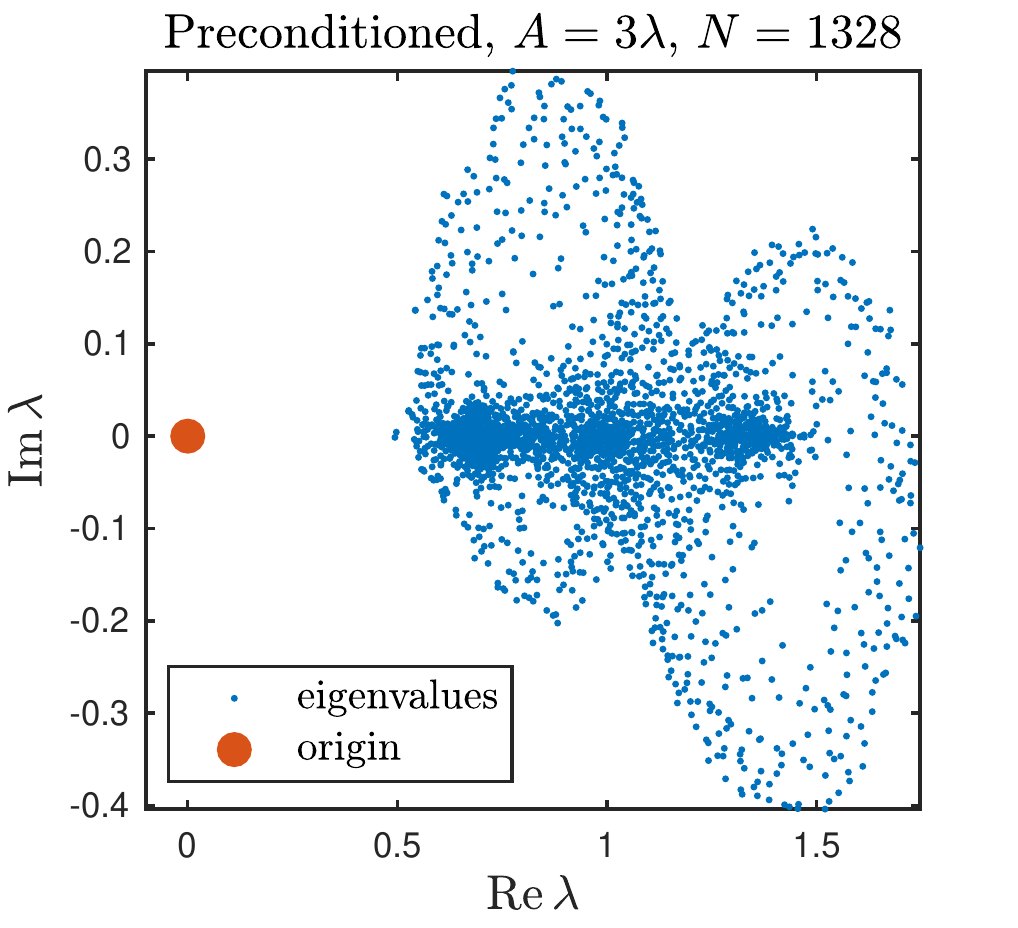}\\
    \includegraphics[width=4.35cm]{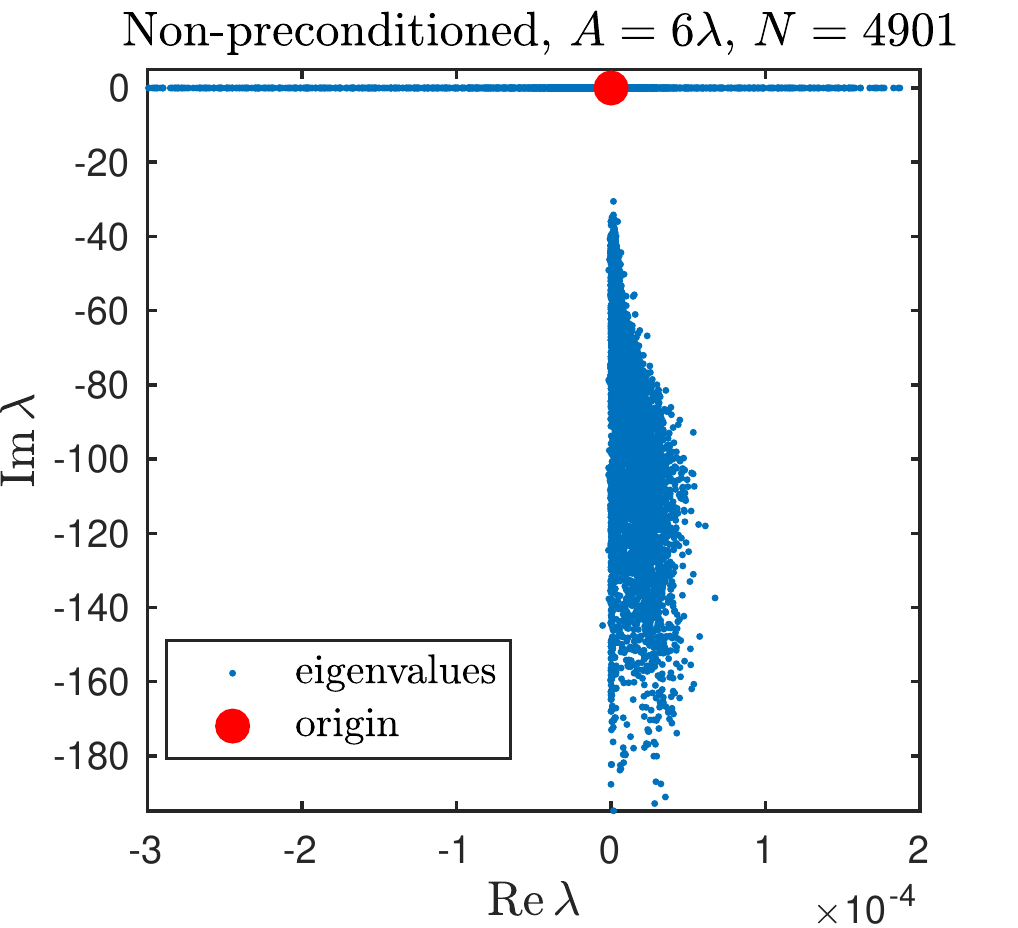}
        \includegraphics[width=4.35cm]{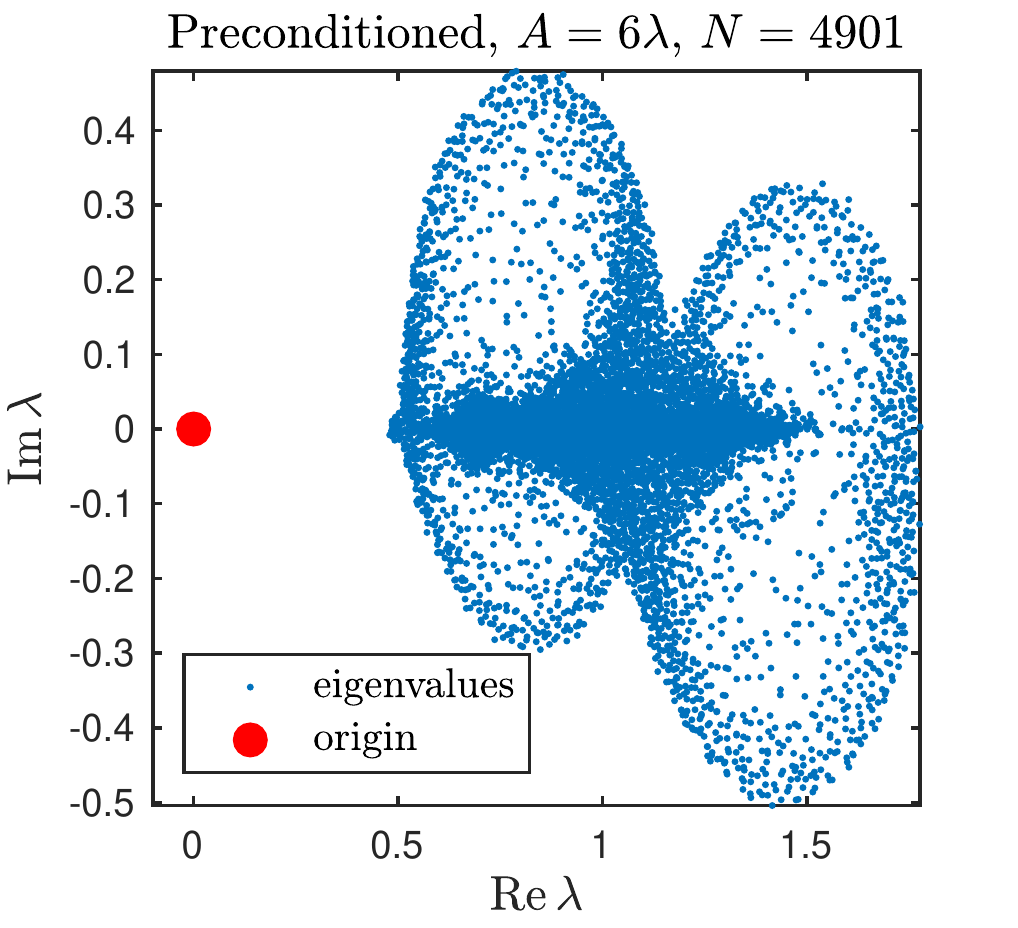}
    \caption{Eigenvalues of the WGF-MoM matrices~\eqref{eq:sys_matrix} (left) and corresponding diagonally preconditioned matrices (right) for the solution of the Sommerfeld half-space problem of Sec.~\ref{sec:SomProbblem} using $A=3\lambda$ (top) and $A=6\lambda$ (bottom). }
    \label{fig:eigen}
    \end{figure}

Interestingly, as in the PEC hemispherical bump example and in all the examples presented in this work, the Jacobi diagonal preconditioner significantly reduces  the number of GMRES iterations required to approximately solve the resulting WGF-MoM algebraic linear system~\eqref{eq:MoMLS}. To examine this fact in more detail we present Fig.~\ref{fig:eigen} which shows the eigenvalues of the non-preconditioned and the Jacobi-preconditioned matrices corresponding to the Sommerfeld dipole problem using window sizes $A=3\lambda$ and $A=6\lambda$ and meshes of similar size ($h\approx 0.3$). Clearly, and unlike the spectra of the non-preconditioned matrices (left), the spectra of the preconditioned matrices (right) are tightly bounded away from the origin, which explains the excellent performance of the preconditioner which reduces from 72 to 13 (resp. 78 to 16) the number of GMRES iterations to achieve the tolerance $10^{-5}$ in the  $A=3\lambda$   (resp.~$A=6\lambda$) case.

Finally, we mention that the LGF code~\cite{panasyuk2009new} took $\sim$8.2 min to carry out the 23320 LGF evaluations needed to produce the reference field~$\elf^{\rm ref}$, while our  (non-accelerated) FORTRAN implementation of the WGF-MoM took just $\sim$5.8 min to produce~$\widetilde\elf^{\rm tot}$ with a relative error~\eqref{eq:rel_error} smaller than 5\%. The calculations were performed on the same machine using one thread in both cases. The runtime difference is more significant when larger numbers of source/target points are considered. Using for instance 9,378 target points (on the surface of the cylinder) and the same~10 sources, the LGF runtime was $\sim$32 min versus just $\sim$7 min of the WGF-MoM. A performance comparison of the WGF-MoM against the  more efficient LGF Strata library~\cite{sharma2021strata} is presented in next section for the solution of a PEC scattering problem. Unfortunately, Strata does not directly produce the entire dyadic LGF, so it could not be used in the examples presented in this section.

 \subsection{PEC obstacle above a dielectric half-space\label{sec:DCIM_comparison}} 

In order to assess the efficiency of the proposed WGF-MoM in this section we compare  its performance against a LGF-MoM for the solution of a planewave scattering problem by a PEC object embedded in the two-layer medium considered in Sec.~\ref{sec:SomProbblem} above. Letting $\mathcal S^S_1$ and $\mathcal D^S_1$ denote the off-surface operators~\eqref{eq:EFIE_pot} and \eqref{eq:MFIE_pot} where integration is performed over the PEC surface~$S$, our WGF formulation of the problem is derived from the field representation 
\begin{align*}
\elf_j(\ner) :=&~k_j^2(\mathcal{S}_{j}\bol v)(\ner)+i\omega\mu_j\{(\mathcal{D}_{j}\bol u)(\ner)+\delta_{1j}(\mathcal{D}^S_{j}\bol w)(\ner)\}\\
\mgf_j(\ner) :=&~k_j^2\{(\mathcal{S}_{j}\bol u)(\ner)+\delta_{1j}(\mathcal{S}^S_{j}\bol w)(\ner)\}-i\omega\epsilon_j(\mathcal{D}_{j}\bol v)(\ner)
\end{align*}
for  $\ner\in\Omega_j$, $j=1,2$, where $\delta_{11}=1$ and  $\delta_{12}=0$, of the EM field scattered by $S$. This representation leads to a second-kind SIE system consisting of a MFIE block that enforces the PEC boundary condition $\hat{\mathbf{n}} \times \mathbf{E}_1=-\hat{\mathbf{n}} \times \mathbf{E}^{\mathrm{src}}$ on~$S$, that is coupled to a windowed M\"uller block that accounts for the transmission condition at the planar dielectric interface $\Gamma$. 

\begin{figure}[h!] 
    \centering
    \includegraphics[width=8.5cm]{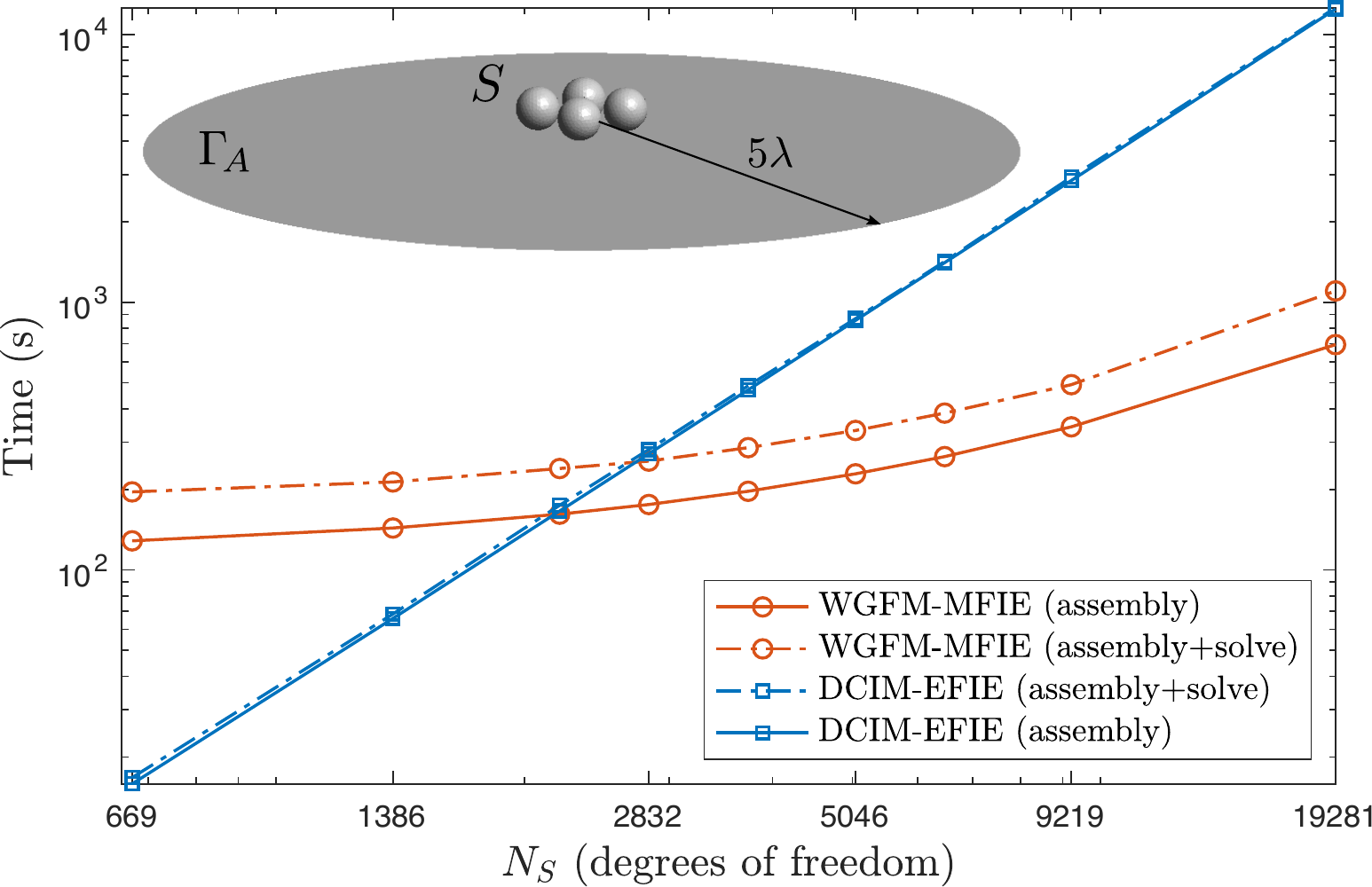}
    \caption{WGF- and LGF-MoM full matrix assembly and linear system solution times in a planewave scattering problem by a PEC obstacle ($S$) above a dielectric half-space for various obstacle discretizations ($N_S$). The LGF-MoM is based on the mixed-potential EFIE formulation~\cite{michalski1990electromagnetic,michalski1990electromagneticII} with the LGF potentials evaluated by means of the Strata library~\cite{sharma2021strata} through its DCIM option.}
    \label{fig:times}
\end{figure}

The LGF-MoM, on the other hand, is based on the mixed-potential EFIE formulation put forth in~\cite{michalski1990electromagnetic,michalski1990electromagneticII}, whereby the scattered field is expressed as $\bold E(\ner) = \int_{S} \bold G(\ner,\ner')\bold J(\ner')\de s'$ everywhere in $\R^3\setminus S$ in terms of the LGF. The singularity cancellation technique~\cite{Sauter2010} is used to treat the singular (static) part of the resulting dyadic and scalar EFIE  kernels. The remaining smooth (dyadic and scalar) parts, which are evaluated using the highly efficient DCIM option in Strata~\cite{sharma2021strata}, are numerically integrated using the standard three-node quadrature rule. A small-size scatterer $S$ placed well above the interface is herein used so as to maintain the accuracy and robustness of the DCIM, which is initialized only once during the whole EFIE matrix assembly. The actual surface used consists of four well separated spheres of equal radius $\lambda/4$ centered at $(\lambda,\pm\lambda,\lambda)/2$ and $(\pm\lambda,\lambda,\lambda)/2$ (see inset in Fig.~\ref{fig:times}).

Figure~\ref{fig:times} displays the time (in seconds) required by each method to assemble the corresponding full system matrices and solve the linear systems for various mesh refinements of~$S$, which are characterized by the number of mesh edges $N_S$. The WGF-MoM system is solved iteratively by diagonally preconditioned GMRES with a tolerance of $10^{-4}$, which converged in less than 27 iterations in all the examples considered.  The LGF-MoM system is in turn solved by direct LU factorization (no speed-up is gained in this case employing GMRES). The effect of the lower dielectric half-space is accounted for with a precision of about 7\% in both cases, which is achieved by suitably selecting the relevant DCIM and WGF parameters. In the latter case the value $A=5\lambda$ is used in the definition of $\Gamma_A$ which is discretized using a mesh of size $h=6.5\times 10^{-2}\lambda$ consisting of 7,238 edges.  The WGF-MoM matrices considered in this comparison have then dimensions $(N_S+\text{14,476})\times(N_S+\text{14,476})$ while the LGF-MoM matrices are much smaller, of dimensions $N_S\times N_S$. Both WGF- and LGF-MoM codes were ran in the same computer using no parallelization of any kind. For the sake of fairness, both codes were written in FORTRAN 90 so that they could share the largest possible number of source code lines. A wrapper of Strata was developed in order to make that possible. 

These results show that, even in this simple setting where no dielectric interface perturbations are considered, the proposed  WGF-MoM outperforms the LGF-MoM for moderately refined (or large) surfaces $S$. This difference in performance between the two methods is mainly explained by the fact that one free-space Green function evaluation costs significantly less than one LGF evaluation, even when the efficient DCIM is employed (that difference is orders of magnitude larger when the LGF is directly evaluated via numerical integration techniques). It so much so that the LGF evaluation cost quickly rises above the cost associated with enforcing the transmission conditions on $\Gamma_{\!A}$ in the WGF approach, which dominates for coarsely refined (or small)  surfaces $S$. 

An even better relative performance of the WGF-MoM is expected when considering, for instance,  structures having localized surface perturbations, like the ones considered in Secs.~\ref{sec:cavity} and~\ref{sec:metalens}, which are particularly cumbersome to deal with  LGF-based methods. In this case the enforcement of the continuity of the total tangential electric and magnetic fields at dielectric interfaces entails evaluation not only of the dyadic LGF itself but also of its curl, further affecting the overall performance of LGF-MoMs. Moreover, in such cases the evaluation of the LGF is hindered by the lack of exponential decay of the spectral LGF when both source and target points lie on the interface between two layers. 
Despite the above mentioned drawbacks, LGF-based SIEs possess remarkable advantages over the WGF approach in certain cases. For instance, structures having several flat layers can be easily handled by the LGF at almost no additional cost, whereas the WGF approach requires the use of additional  surface currents at each one of the interfaces, even when they do not contain any perturbation. Similarly, problems involving moderate numbers of small-area inclusions entail few LGF evaluations, and hence can be easily treaded by this  approach.

 \begin{figure}[!t]
    \centering
    \includegraphics[width=9cm]{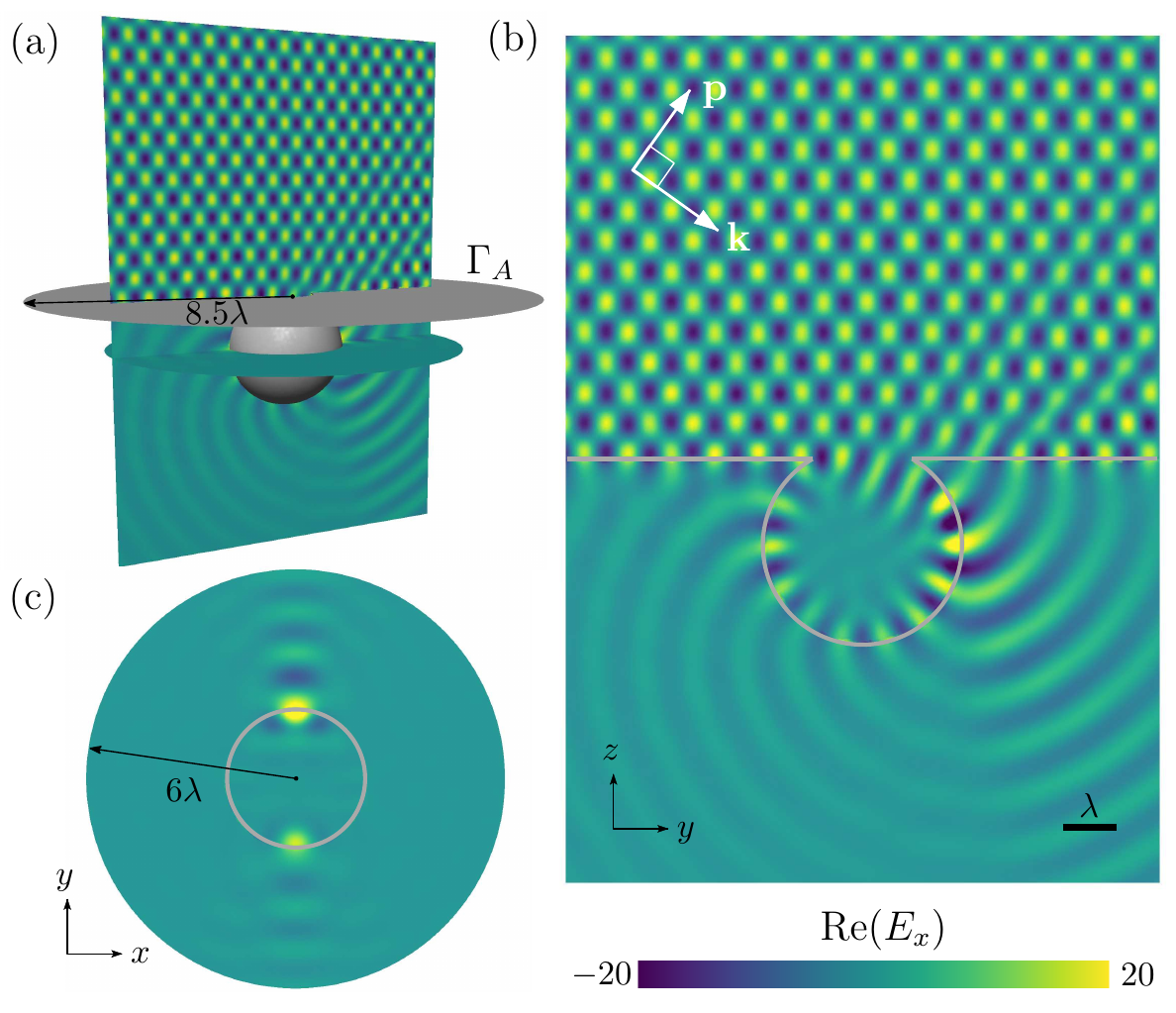}
    \caption{Real part of the $x$-component of the total electric field solution of the problem of scattering of a TE-polarized incident planewave impinging on a spherical cavity in a dielectric half-space, produced by the proposed WGF-MoM with $A=8.5\lambda$. Truncated surface (grey) and field view on the $yz$-plane (a),  $yz$-plane (b), and  $xy$-plane (c) within the region $\{w_A(\ner)=1\}$.}
    \label{fig:cavity}
\end{figure}

 \subsection{Cavity in a dielectric half-space\label{sec:cavity}} 
Next, we consider the problem of scattering of a TE-polarized planewave~\eqref{eq:PW} that impinges from $\Omega_1$ on a large spherical-shaped cavity in a dielectric half-space (see Fig.~\ref{fig:cavity}). Once again we use dimensionless physical parameters, with $k_1=1.4$, $k_2=1$, and $\lambda=2\pi/k_2$. The radius of the spherical sector embedded in the lower half-space is $2\lambda$, giving rise to a disk aperture of radius $\lambda$. A grazing angle of $\pi/5$ was used in this example, so that, in absence of the spherical cavity, total internal reflection would have taken place at the planar interface. The truncated locally-perturbed surface $\Gamma_A$, as well as the real part of the of $x$-component of the total electric field ($E_x$), are shown in Fig.~\ref{fig:cavity}. Figures~\ref{fig:cavity}(b) and \ref{fig:cavity}(c) display the real part of electric field on the $yz$- and $xy$-planes, respectively, within the region $\{w_A(\ner)=1\}$ where formula~\eqref{eq:windEMfld} (used to the produce the fields) yields an accurate field evaluations. The approximately uniform triangular mesh of $\Gamma_A$ used in the MoM-produced fields displayed in Fig.~\ref{fig:cavity},  corresponds to $A=8.5\lambda$, $h=0.15\lambda$, and consists of a total of $N=97321$ edges. The linear system~\eqref{eq:MoMLS} was solved by a Jacobi-preconditioned GMRES solver which converged in 50~iterations to the prescribed $10^{-4}$ tolerance. Note that the field $E_x$ plotted in Fig.~\ref{fig:cavity}(b) looks continuous across the material interface, as it should be, and that transmission to the lower half-space takes place only within the cavity, due to the total internal reflection incidence.

\begin{figure}[h!] 
    \centering
    \includegraphics[width=9cm]{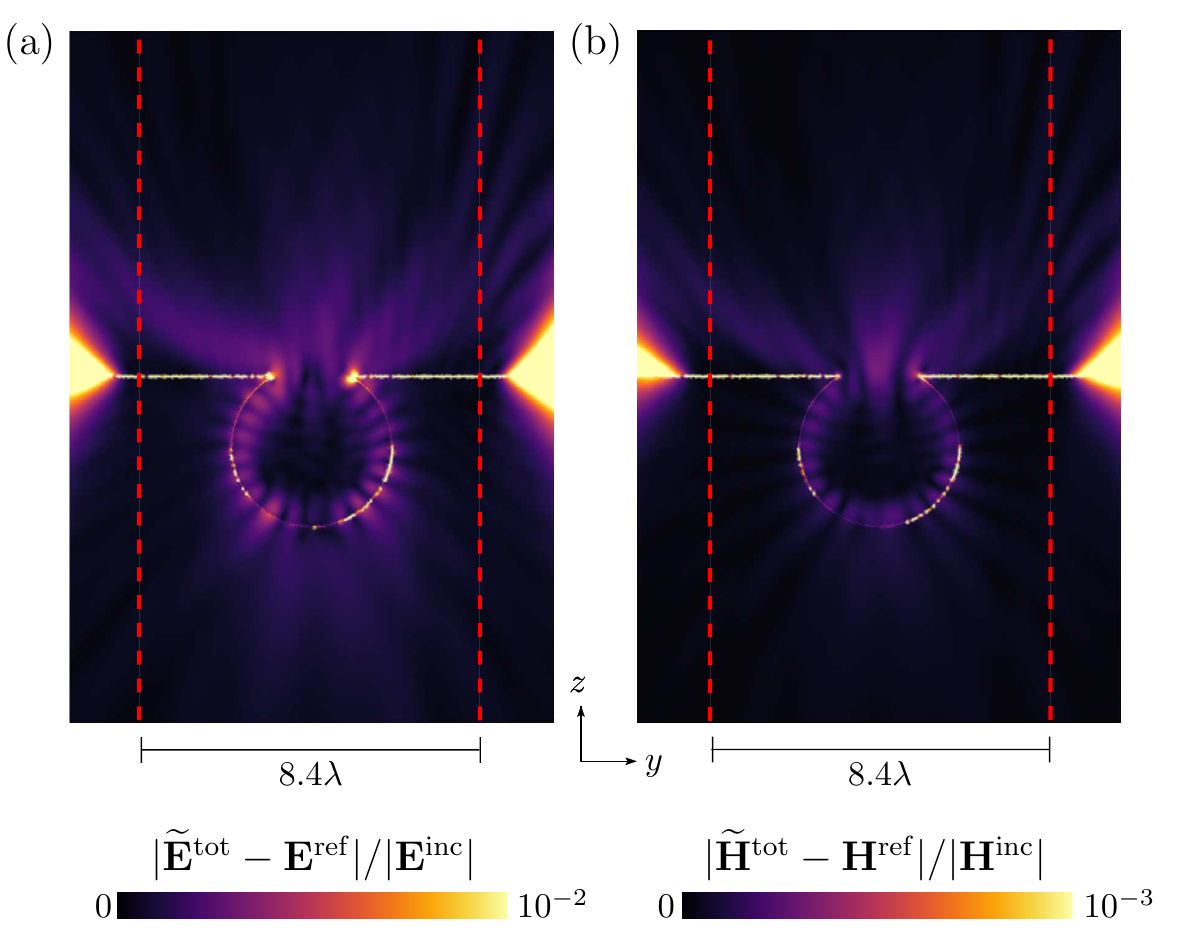}
    \caption{Absolute errors in the total electric (a) and magnetic (b) fields, normalized by the corresponding (constant) incident-field amplitude, in the WGF-MoM solution of the cavity problem using $A=6\lambda$. The reference EM field $(\elf^{\rm ref}, \mgf^{\rm ref})$ was produced using a larger window size of $A=8.5\lambda$ and a similar mesh size. The dashed red lines mark the boundary of $\{\ner\in\R^3:w_{6\lambda}(\ner)=1\}$.}
    \label{fig:cavityError}
\end{figure}

In order to assess the field errors in the solution of the cavity problem, we solve the same cavity problem using a slightly smaller window size, $A=6\lambda$, and compare it to the solution shown in Fig.~\ref{fig:cavity} corresponding to $A=8.5\lambda$. Both surface meshes have approximately the same size $h=0.15\lambda$ and an almost identical number of GMRES iterations were needed to achieve the $10^{-4}$ tolerance. The (normalized) absolute errors in the total electric and magnetic fields are displayed in Figs.~\ref{fig:cavityError}(a) and~(b), respectively, in a portion of the $xy$-plane. The boundary of the strip $[-6\lambda,6\lambda]\times \R$ contained in $\{w_A(\ner)=1\}$ is marked by the red vertical dashed lines. The severe loss of accuracy occurring outside $\{w_A(\ner)=1\}$ is clearly observed in those figures, specially near the interface~$\Gamma$. In particular, this example shows that the WGF-MoM does not directly produce correct far fields. As mentioned in Sec.~\ref{sec:window}, a simple remedy to this problem is to map the correct near fields, produced by the WGF-MoM, to the far field. This can be done by means of Stratton-Chu formula based on the LGF integrating over a surface enclosing the perturbation, and then replacing the kernels by their leading-term asymptotic expansion as $|\ner|\to\infty$.

\begin{figure}[!t]
  \centering
  \includegraphics[width=9cm]{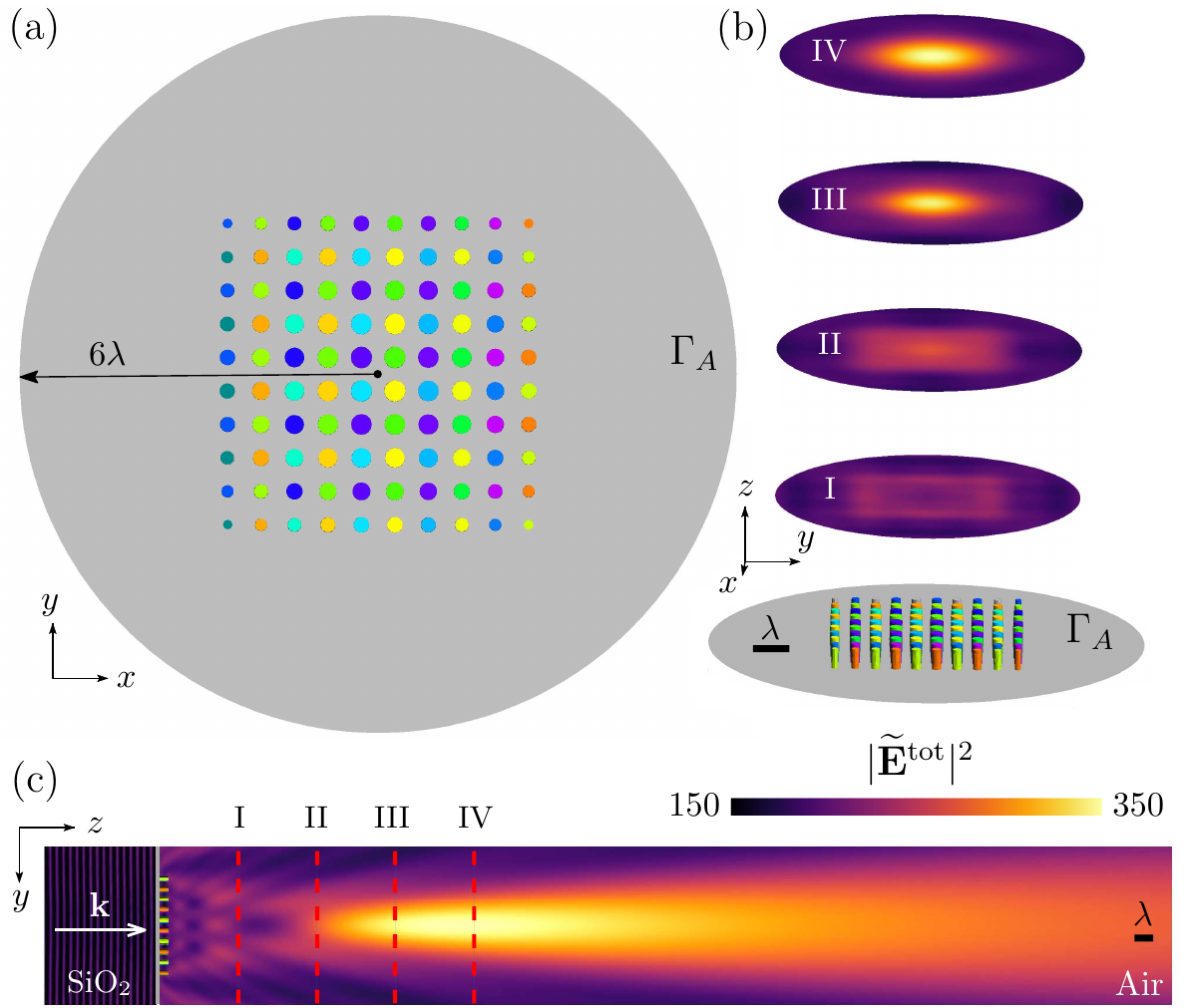}
  \caption{Refraction of a normally incident  plane EM wave by an all-silica metasurface. (a) Top view of the metasurface consisting of a $10\times 10$ array of sub-wavelength nano-rods of variable radii. (b)-(c) Two different views of the total electric field intensity $|\widetilde\elf^{\rm tot}|^2$ obtained by means of the proposed WGF-MoM, demonstrating the light focusing effect of the metasurface. The dashed red lines in~(c) mark the location of the surfaces where $|\widetilde\elf^{\rm tot}|^2$ is plotted in~(b).}
  \label{fig:metasurface}
  \end{figure}

\subsection{All-silica metasurface}\label{sec:metalens}
Our next example presents the full 3D solution of a problem of scattering by an all-silica (${\rm SiO}_2$) metasurface (cf.~\cite{li2017free}) consisting of an array of $10\times10$ nano-rods of sub-wavelength radii, ranging from $36.7
$ nm to $77.5$ nm, and a fixed height of 364 nm, which is illuminated from  $\Omega_2$ by a normally incident plane EM~wave~\eqref{eq:PW} with $\bold p=(1,1,1)$ and $\bold k =k_2(0,0,1)$. Figure~\ref{fig:metasurface}(a) displays the surface $\Gamma_A$  utilized in this example, where the window size $A=6\lambda=3.9~\mu{\rm m}$ ($\lambda=650~{\rm nm}$) is used. (The nano-rods are shown in various colors for visualization purposes.) The nano-rod radii follow a parabolic profile that effectively steers the direction of the transmitted light by controlling the phase change as it penetrates $\Gamma_A$, thus achieving the focusing effect demonstrated in Figs.~\ref{fig:metasurface}(b) and~(c), that display the electric field intensity at surfaces parallel to the $yz$- and $xy$-plane, respectively. The surface mesh employed in this example, which was properly refined so as to account for the numerous small-scale features, comprised a total of 17,822 nodes and 53,335 edges. The corresponding SIE solution was produced using a Jacobi-preconditioned GMRES solver, which converged  in 45 iterations to a tolerance of $10^{-4}$. Our (unaccelerated) FORTRAN OpenMP-parallelized WGF-MoM implementation took $\sim$6 min in constructing the system matrix, $\sim$15 min in solving the linear system, and  $\sim$35 min in producing the high-fidelity fields shown in Fig.~\ref{fig:metasurface} using a workstation with 48 cores (96 threads, dual Xeon Gold 6240R) and 500 GB of memory.

\subsection{Plasmonic solar cell}\label{sec:plasmon}
\begin{figure}[!t]
  \centering
  \includegraphics[width=9cm]{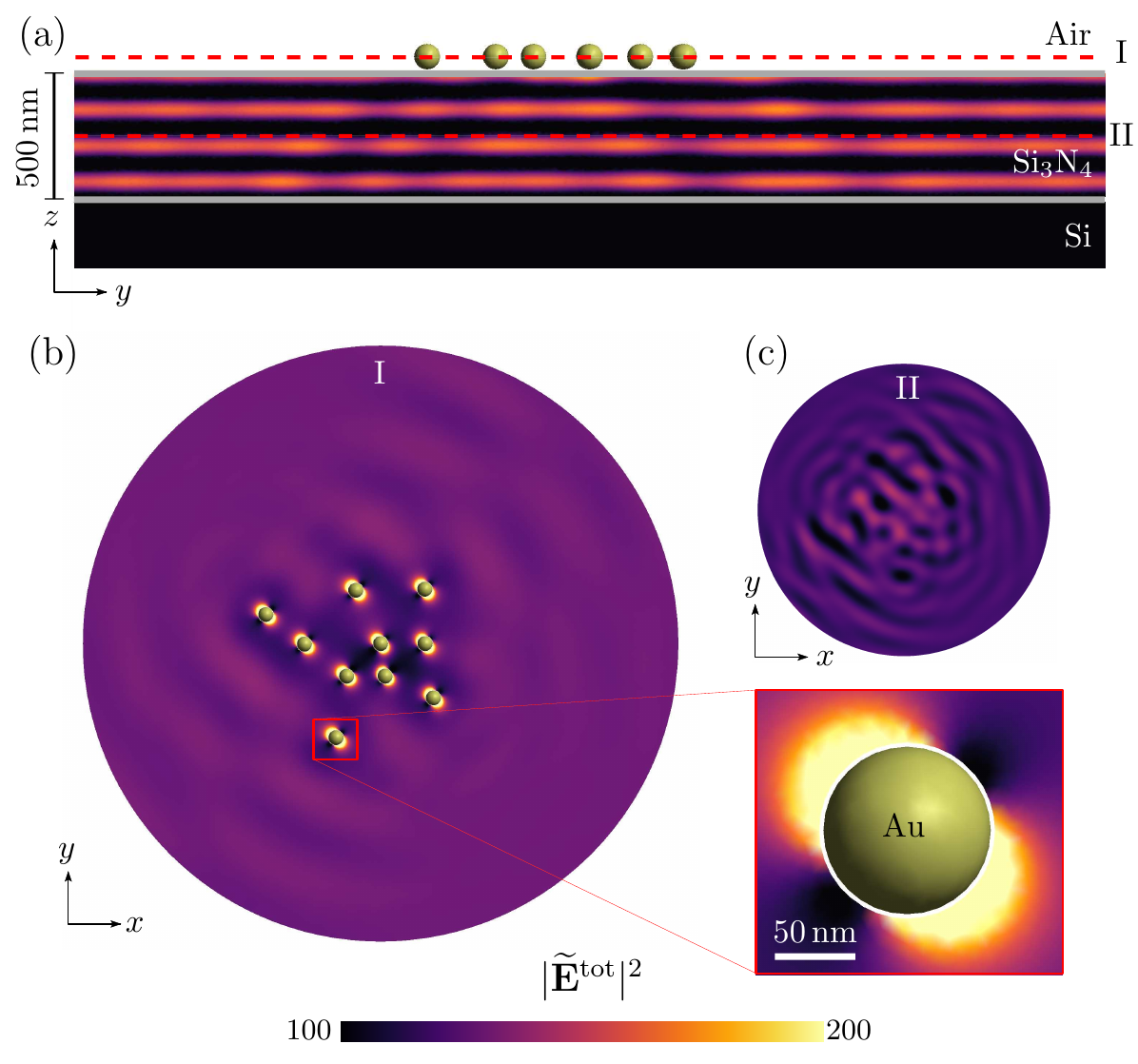}
  \caption{Planewave illumination of a plasmonic solar cell structure consisting of 10 gold nanoparticles randomly placed on top of a silicon nitride film backed by a silicon substrate. (a) Electric field intensity $|\widetilde \elf^{\rm tot}|^2$ cross-section plot obtained by means of the proposed WGF-MoM. (b)-(c) Plots of $|\widetilde\elf^{\rm tot}|^2$ on the planes marked by the dashed red lines in~(a). The plasmon resonance field enhancement around a gold nanoparticle is clearly visible in the zoomed inset figure in~(b).}
  \label{fig:solarCell}
  \end{figure}

In the last example of this paper  we consider a plasmonic solar cell structure~\cite{atwater2011plasmonics} consisting of~10 gold nanoparticles of diameter $100$ nm lying on top of a 500 nm thick silicon nitride (${\rm Si_3N_4}$) film backed by a silicon (Si) substrate (see Fig.~\ref{fig:solarCell}(a)). The structure is illuminated by a normally incident planewave~\eqref{eq:PW} coming from above at $\lambda = 572$ nm and polarized according to $\bold p = (1,1,1)$. The total electric field intensity~$|\widetilde\elf^{\rm tot}|^2$ produced by the proposed WGF-MoM is shown in Fig.~\ref{fig:solarCell}. The frequency used in this example excites plasmon resonances in the metallic nanoparticles leading to a strong local field enhancement around the metallic nanoparticles~\cite{maier2005plasmonics}, as can be observed in the zoomed inset figure in Fig.~\ref{fig:solarCell}(b). 

The presence of multiple penetrable interfaces made it necessary to generalize the SIE formulation presented above in Sec.~\ref{sec:SIEFormulation}. In detail, letting $\Omega_1$, $\Omega_2$, $\Omega_3$ and~$\Omega_4$ denote  the subdomains occupied by air ($\epsilon_1=\epsilon_0$), silicon nitride ($\epsilon_2=2.0483\epsilon_0$), silicon ($\epsilon_3=(4.0191+0.031373i)\epsilon_0$), and the gold nanoparticles ($\epsilon_4=(0.33221+2.74i)\epsilon_0$), respectively, we express the EM fields $(\elf_j,\mgf_j)$ in $\Omega_j$ as in~\eqref{eq:SC} but in terms of the off-surface operators $\mathcal S_j$~\eqref{eq:EFIE_pot} and $\mathcal D_j$~\eqref{eq:MFIE_pot} defined by integrals over $\p\Omega_j$, $j=1,\ldots,4$. Enforcing then the continuity of the total tangential fields $\bold n\times (\elf_j+\elf^{\rm src},\mgf_j+\mgf^{\rm src})$ at each of the interfaces $\p\Omega_j$, using as source field $(\elf^{\rm src},\mgf^{\rm src})$  the total EM field solution of the problem of scattering by the three-layer  structure (without the nanoparticles) (see~\cite[Sec.~2.1.3]{Chew1995waves}), we arrive at a $3\times 3$ block second-kind SIE system that is windowed and discretized using the MoM presented in Sec.~\ref{sec:MoM}. The approximate total fields~$(\widetilde\elf^{\rm tot},\widetilde\mgf^{\rm tot})=(\widetilde\elf_j+\elf^{\rm src},\widetilde\mgf_j+\mgf^{\rm src})$, $j=1,\ldots,4$, are retrieved by windowing the corresponding field representation formulae.

The two planar triangular meshes used in this example comprise 35,785 edges each while the total number of edges in the spherical meshes amounted to 11,337. The planar meshes were suitably refined near the bottom tip of the spheres to properly account for possible nearly singular integration issues. The linear system was solved by means of GMRES, which converged in 28 iterations to the prescribed tolerance ($10^{-4}$). The overall  time needed by our OpenMP-parallelized WGF-MoM implementation to produce the three plots of $|\widetilde\elf^{\rm tot}|^2$ presented in Fig.~\ref{fig:solarCell}, on the horizontal planes at $z=50~\rm{nm}$ and $z=-250~{\rm nm}$ in~(b) and~(c), respectively, as well as on the vertical plane $\{x=0\}$ in~(a), was around 64 min on the aforementioned 48-core machine. 

\section{Conclusions and Discussion}

This paper presents a  SIE method for EM scattering by locally perturbed planar layered media. The proposed methodology, which extends the WGF method put forth in~\cite{Bruno2015windowed,bruno2017windowed,perez2017windowed} for the (scalar) Helmholtz equation, does not entail evaluation of any  Sommerfeld integrals thus avoiding their inherent costs and challenges that they pose, but at the expense of adding new unknowns on the planar interfaces and requiring a larger linear system that must be solved. The method leverages an indirect second-kind M\"uller SIE formulation featuring weakly-singular integral operators expressed in terms of free-space  Green functions. Upon windowing the integral kernels and applying a standard Galerkin-MoM discretization based on RWG basis functions, well-conditioned linear systems amenable to be solved iteratively by GMRES, are obtained. The resulting methodology exhibits both second-order convergence (in the near fields) as the mesh size in decreased, and high-order convergence (super-algebraic) as the window size is increased, as demonstrated by a thorough set of comparative examples. A number of challenging problems including scattering by cavities, metasurfaces, and plasmonic solar cell structures, further validate and showcase the capabilities of the proposed WGF-MoM.  
It is worth mentioning, however, that larger-scale and more realistic metasurface and solar cells configurations than the ones considered here, inevitably require use of fast algorithms such as the fast multipole method~\cite{song1997multilevel} or $\mathcal H$-matrices~\cite{hackbusch2015hierarchical}. 

This work certainly opens up a number of possible follow-up research directions. For example, an accurate SIE solver capable of handling more general metasurface designs requires proper handling of multi-material junctions. As in the 2D case~\cite{pestourie2018inverse,jerez2017multitrace}, this could be accomplished within our 3D WGF-MoM framework by employing a second-kind single-trace formulation~\cite{claeys2017second}. The robust low-frequency behavior of M\"uller's formulation reported in~\cite{yla2005well}, on the other hand, brings about the idea of extending the proposed methodology to the time domain by suitably combining it with convolution quadrature schemes, as was done in the 2D case in~\cite{labarca2019CQ}. 


\appendices
\section{Hemispherical bump problem: Mie series solution}\label{app:A}
We here make use of the classical Mie series solution and the theory of images, to produce the exact solution of the problem of scattering of a plane EM wave by a hemispherical PEC bump on top of a PEC half-space. Consider an incident plane EM wave, with grazing angle $\alpha$, given by 
\begin{equation*}
	\mathbf{E}_\alpha^{\inc}(\ner ) = E_{0}\e^{ik\left(y\cos \alpha-z\sin \alpha\right)}
	\begin{cases}
	\hat\nex \qquad\,\,\text{if TE polarized,}\\
	\hat\ney\sin \alpha+\hat\nez\cos \alpha \\
	\,\qquad\quad\text{if TM polarized.}\\
	\end{cases}
\end{equation*}
Using the theory of images we have that the resulting total electromagnetic field, which satisfies
$\nor \times \elf^{\rm tot}=\bol 0$  on $\Gamma$, can be expressed as 
$$
\elf^{\rm tot}=\elf_\alpha^{\inc}+\elf^{\rm Mie}+\widetilde\elf^{\inc}+\widetilde\elf^{\rm Mie}
$$
where
\begin{equation*}
\widetilde\elf^{\inc} = 
\begin{cases}
-\mathbf{E}_{-\alpha}^{\inc} &\quad\text{if TE polarized}\\
\mathbf{E}_{-\alpha}^{\inc} &\quad\text{if TM polarized}\\
\end{cases}
\end{equation*}
and where $\elf^{\rm Mie}$ and $\widetilde\elf^{\rm Mie}$ are the well-known Mie series solution of the problem of scattering of an entire PEC sphere, with the same radius as the bump, by $\elf^\inc$ and $\widetilde\elf^\inc$, respectively.

\appendices

\section*{Acknowledgment}

The authors gratefully acknowledge the support by FONDECYT (Fondo Nacional de Desarrollo Científico y Tecnológico) Chile, through Grant No. 11181032. Rodrigo Arrieta thanks ANID (Agencia Nacional de Investigación y Desarrollo), Subdirección de Capital Humano/Magíster Nacional/año 2021/folio 22211890, for funding his postgraduate studies.

\bibliographystyle{IEEEtran}
\bibliography{References}{}

\end{document}